\documentclass[11pt]{article}

\usepackage{amsmath,amsthm,verbatim,amssymb,amsfonts,amscd, graphicx, tikz-cd, pinlabel, thm-restate, xcolor}
\usepackage{graphics}
\theoremstyle{plain}
\newtheorem{theorem}{Theorem}

\newtheorem{lemma}{Lemma}

\newtheorem{proposition}{Proposition}

\theoremstyle{definition}

\begin{document}

\title{Representing smooth 4-manifolds as loops in the pants complex}
\author{Gabriel Islambouli and Michael Klug}
\maketitle

\section{Overview}
Simplicial complexes associated to curves on a surface play a central role in 2- and 3-manifold topology, particularly in the study of mapping class groups and Heegaard splittings. Recently, Kirby and Thompson \cite{KT} pushed these techniques into dimension four, assigning a loop in the cut complex to a trisected 4-manifold. The aim of this paper is to, in some sense, reverse this. In particular, given a loop in the pants complex, $L$, we show how to uniquely build a closed smooth 4-manifold $\mathcal{X}^4_C(L)$. Our main theorem is that all such manifolds arise in this fashion.

\begin{restatable*}{theorem}{existence}
\label{thm:existence}

	For every closed, smooth, orientable 4-manifold $X^4$, there exists a closed loop $L$ in $\mathcal{P}(\Sigma)$ so that $X$ is diffeomorphic to $\mathcal{X}_C^4(L)$.

\end{restatable*}

In their proof of the finite presentability of the mapping class group \cite{HT}, Hatcher and Thurston sketch a proof that the pants complex is simply connected. This result was later fully fleshed out in work of Hatcher \cite{AH}. As our main theorem associated a loop to any 4-manifold, it is natural to ask what the disk it bounds represents. Viewing 4-manifolds from this perspective yields a natural proof of the following theorem, originally due to Pontrjagin and Rohlin \cite{VR}, which is our main application.

\begin{restatable*}{theorem}{cobordismGroup}
\label{thm:cobordismGroup}
Every smooth, oriented, closed manifold is cobordant to $\coprod_m \mathbb{C}P^2 \coprod_n \overline{\mathbb{C}P}^2$. 
\end{restatable*}

Our proof here follows along the lines of recent work of Gay \cite{DG2}, in which the author proves the same theorem by associating a loop of smooth functions on a surface to a 4-manifold. The similarity in these arguments suggests that the pants complex of a surface $\Sigma$ is a good discrete model for space of smooth functions on $\Sigma$. Our proof relies on the simple-connectivity of the pants complex which originally was proven using properties of generic smooth functions on surfaces. Nevertheless, there now exist multiple proofs of the simple-connectivity of the pants complex which rely on different techniques \cite{BK} \cite{BP} which give rise to alternative paths to the theorem, some of which (after using \cite{BW}) are quite elementary.

We also use our correspondence to gain insight into the structure of the pants complex. In particular, given a loop $L$ in the pants complex, we define an invariant $\sigma(L)$, which is the signature of the 4-manifold associated to $L$. This may be calculated using information only of the 1-skeleton of the pants complex, but contains information about possible disks that this loop can bound. In particular, we obtain the following proposition.

\begin{restatable*}{proposition}{loopTriangleBound}
Let $L$ be a loop in the pants complex with $\sigma(L) = n$, then any disk bounded by $L$ must contain at least $n$ 3S-triangles.
\end{restatable*}

\section{The pants complex}

We briefly discuss the pants complex of a surface following \cite{AH}. Let $\Sigma$ be a connected closed orientable surface of genus greater than or equal to 2.  A \textbf{pants decomposition} of $\Sigma$ is a set of $3g-3$ simple closed curves on $\Sigma$ such that cutting $\Sigma$ along these curves results in a disjoint union of $2g-2$ 3-punctured spheres (pairs of pants).  Two pants decompositions of $\Sigma$ are considered the same if the curves are isotopic.  We will be considering the 2-complex $\mathcal{P}(\Sigma)$, called the pants complex of $\Sigma$, whose vertices correspond to isotopy classes of pants decompositions of $\Sigma$. If $\Sigma$ is instead a torus, a pants decomposition is just an isotopy class of an essential curve.  

\begin{figure}
    \centering

        \includegraphics[scale=.3]{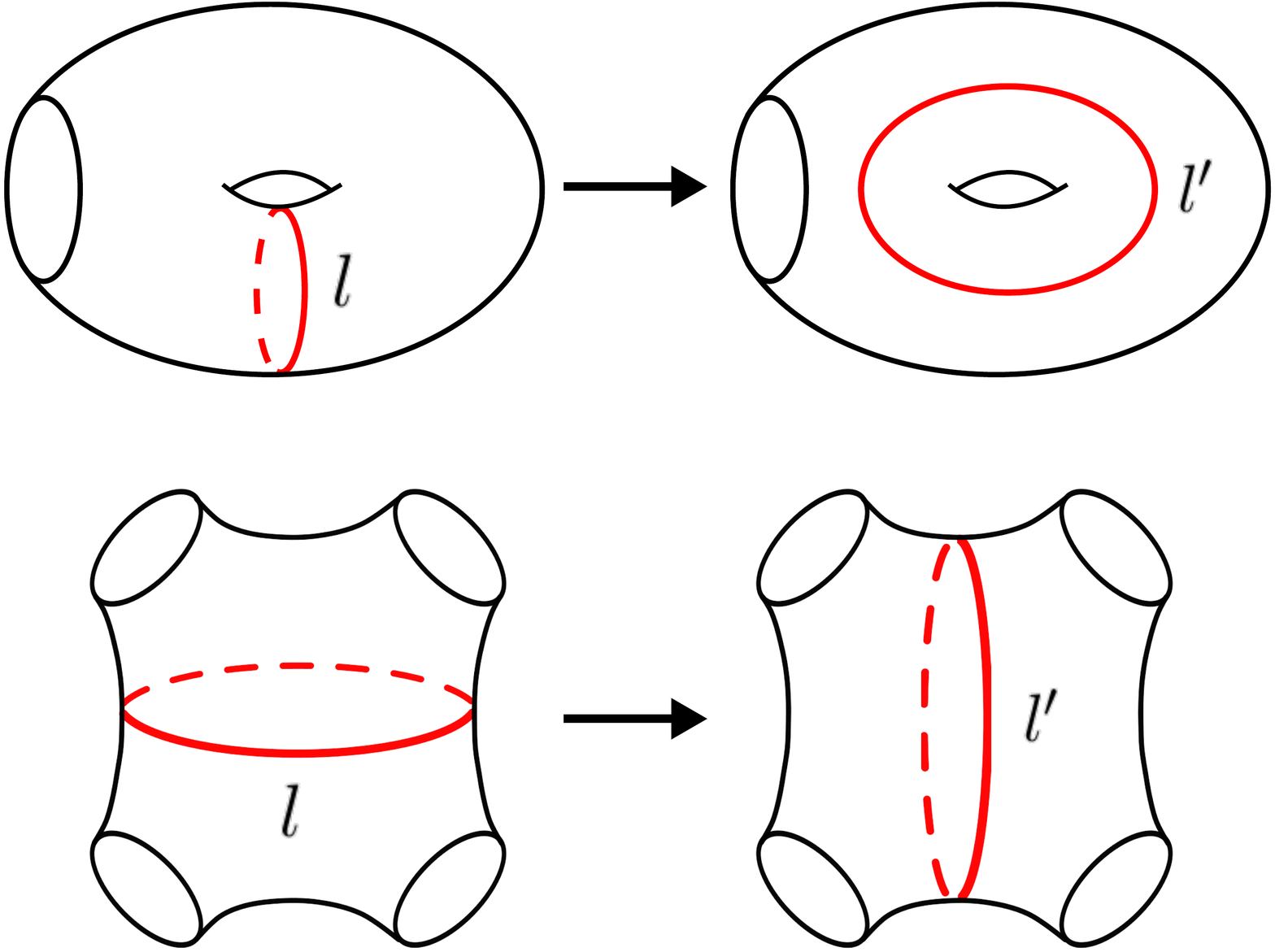}
    \caption{Top: An S-move in the pants complex. Bottom: and A-move in the pants complex}
    \label{fig:PantsComplexMoves}
\end{figure}

There are two types of edges in $\mathcal{P}(\Sigma)$: S-edges (``S" for stabilization) and A-edges (``A" for associative), which can be seen in Figure \ref{fig:PantsComplexMoves}. First note that if one removes neighbourhoods of all but one of the curves in a pants decomposition, one is left with $2g-3$ pants and either a 4-punctured sphere or a once punctured torus. Two pants decompositions $P_1$ and $P_2$ are connected by an $S$-move if all but one of the curves in $P_1$ are the same as in $P_2$ and the curves that differ intersect each other in exactly one point on a once punctured torus component.  Two pants decompositions $P_1$ and $P_2$ are connected by an A-edge if all but one of the curves in $P_1$ are the same as in $P_2$ and the curves that differ intersect each other in exactly two points on a 4-punctured sphere.

The 2-cells in $\mathcal{P}(\Sigma)$ which come in 5 different types and are shown in Figures 2-6. With this choice of 2-cells the following theorem holds.

\begin{theorem}(Hatcher \cite{AH})
\label{thm:simplyConnected}
$\mathcal{P}(\Sigma)$ is connected and simply-connected.  
\end{theorem}

\begin{figure}
    \centering
    \labellist
		\pinlabel {\LARGE{\textbf{$3A$}}} at 593 310
	\endlabellist
        \includegraphics[scale=.3]{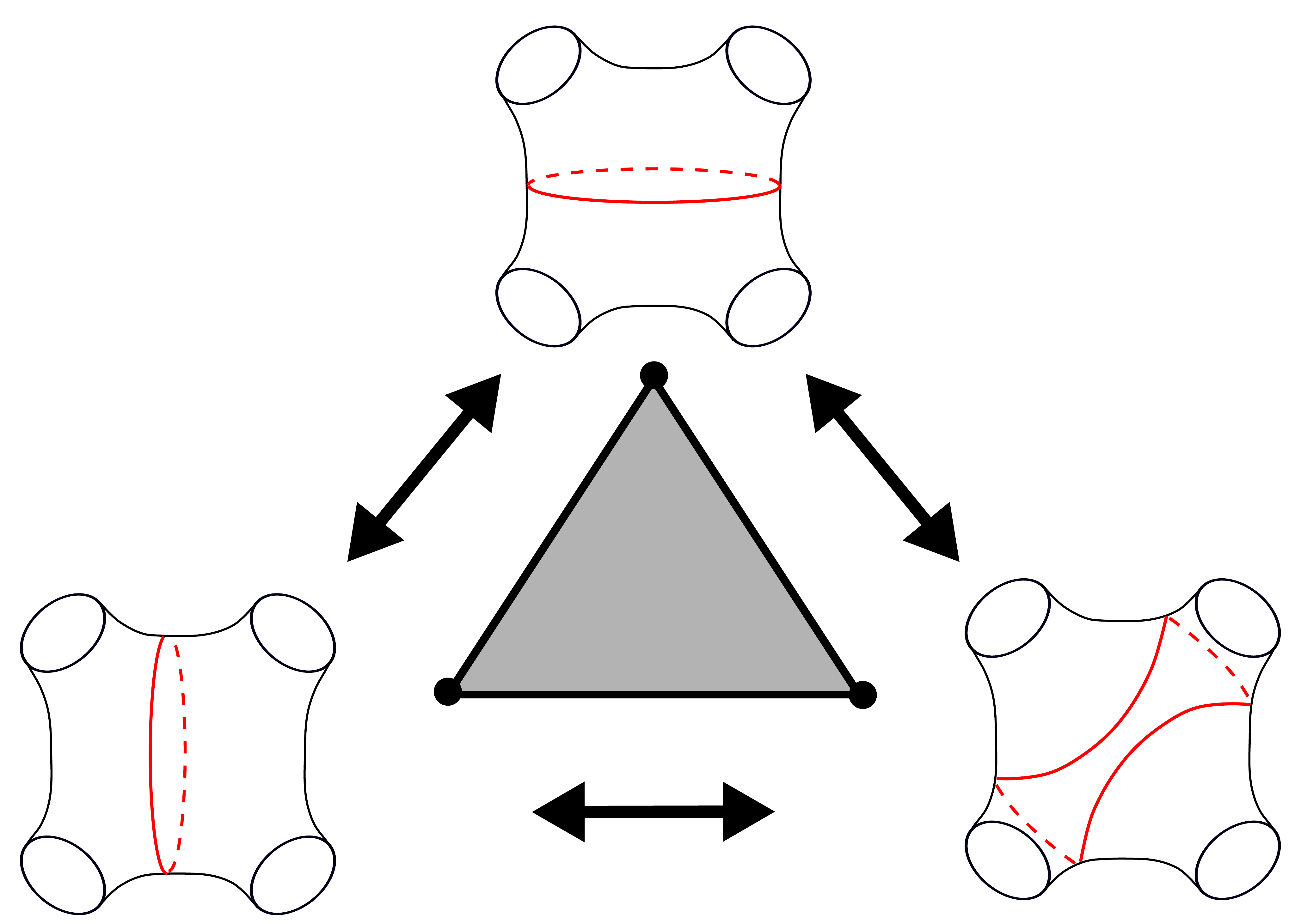}
    \caption{A 3A cycle which is filled in by a triangle.}
    \label{fig:3A}
\end{figure}

\begin{figure}
    \centering
    \labellist
		\pinlabel {\LARGE{\textbf{$3S$}}} at 620 340
	\endlabellist
        \includegraphics[scale=.3]{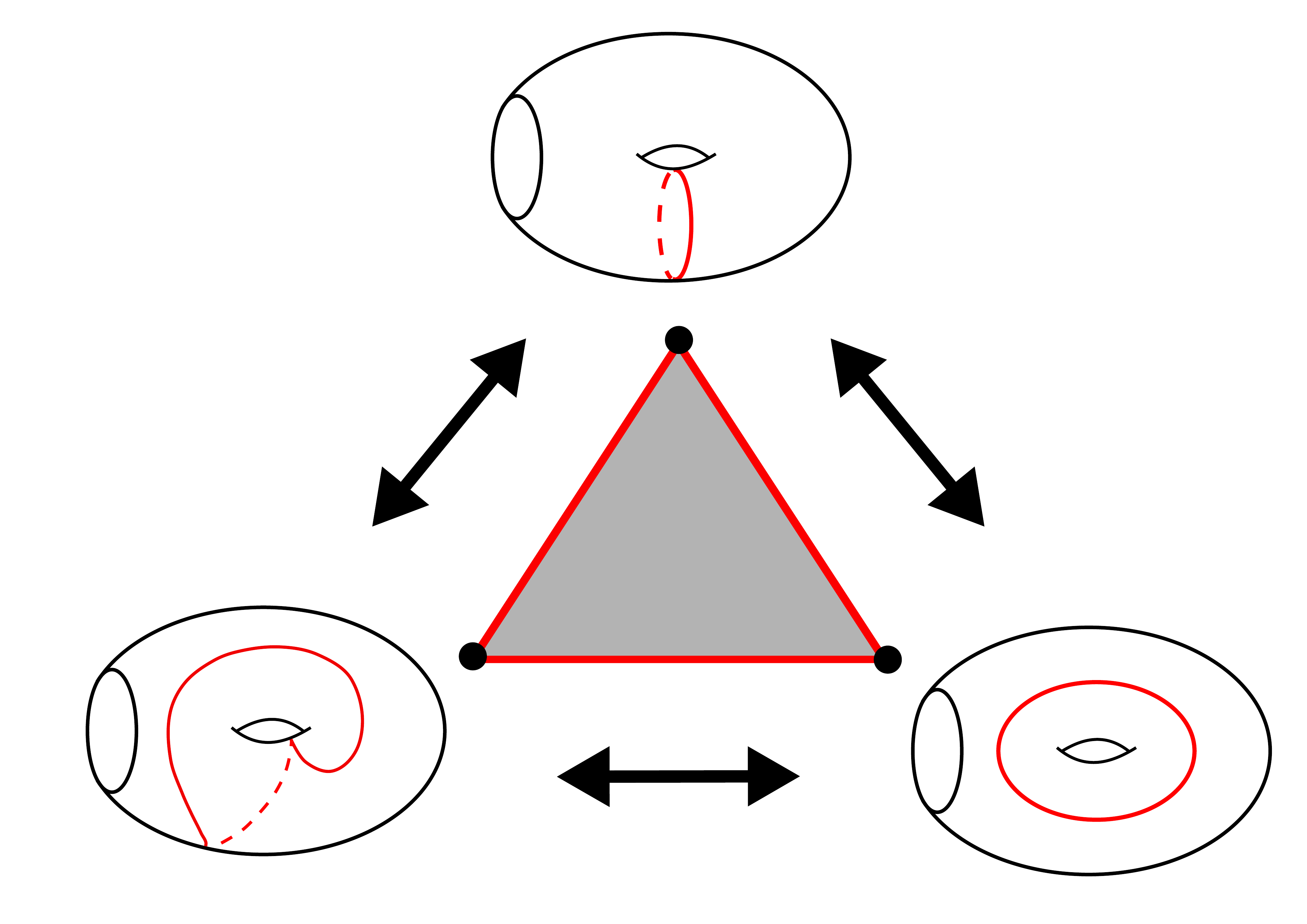}
    \caption{A 3S cycle which is filled in by a triangle.}
    \label{fig:3S}
\end{figure}

\begin{figure}
    \centering
    \labellist
		\pinlabel {\LARGE{\textbf{$5A$}}} at 560 460
	\endlabellist
        \includegraphics[scale=.3]{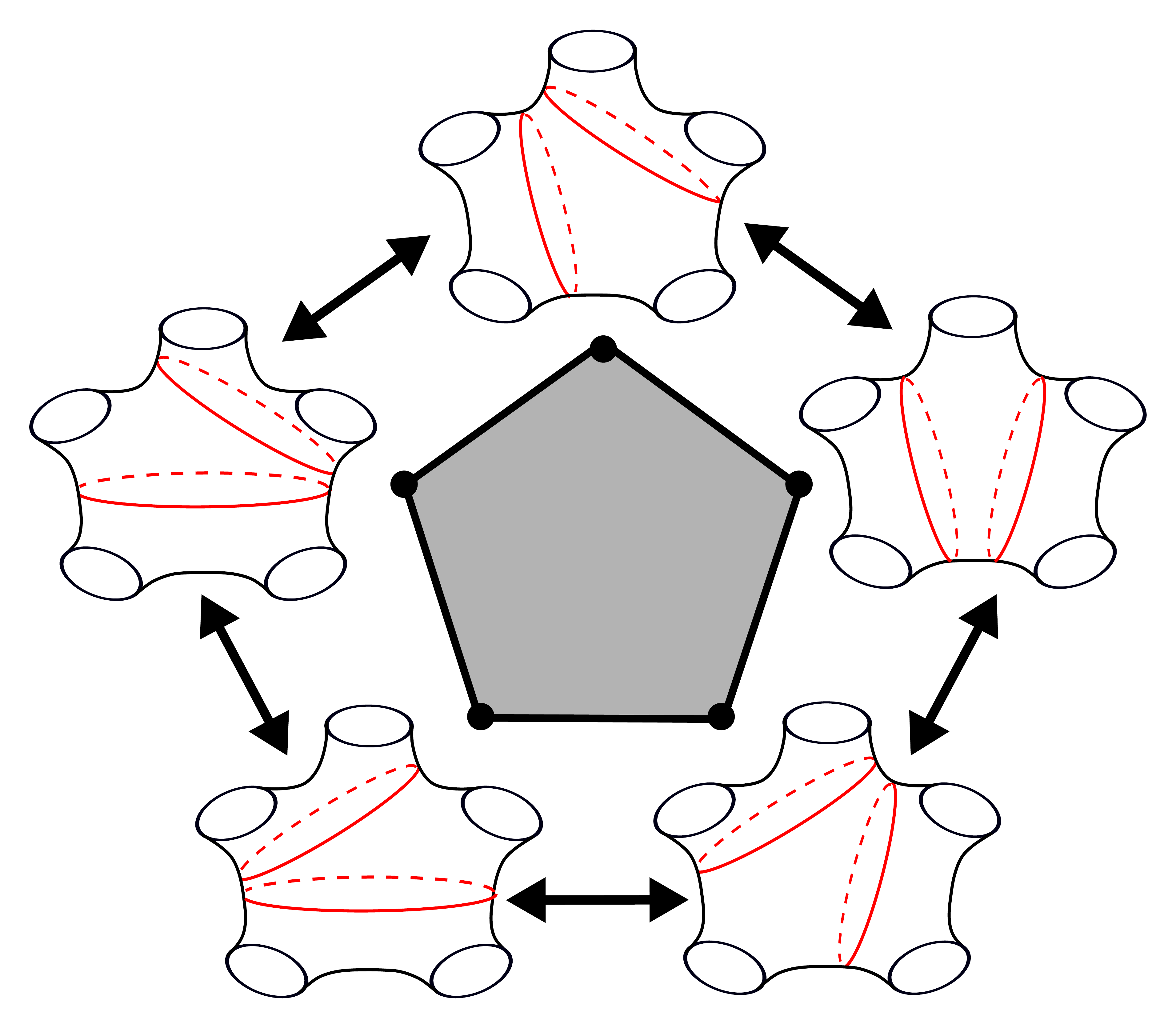}
    \caption{A 5A cycle which is filled in by a pentagon.}
    \label{fig:5A}
\end{figure}

\begin{figure}
    \centering
    \labellist
		\pinlabel {\LARGE{\textbf{$6AS$}}} at 585 440
	\endlabellist
        \includegraphics[scale=.3]{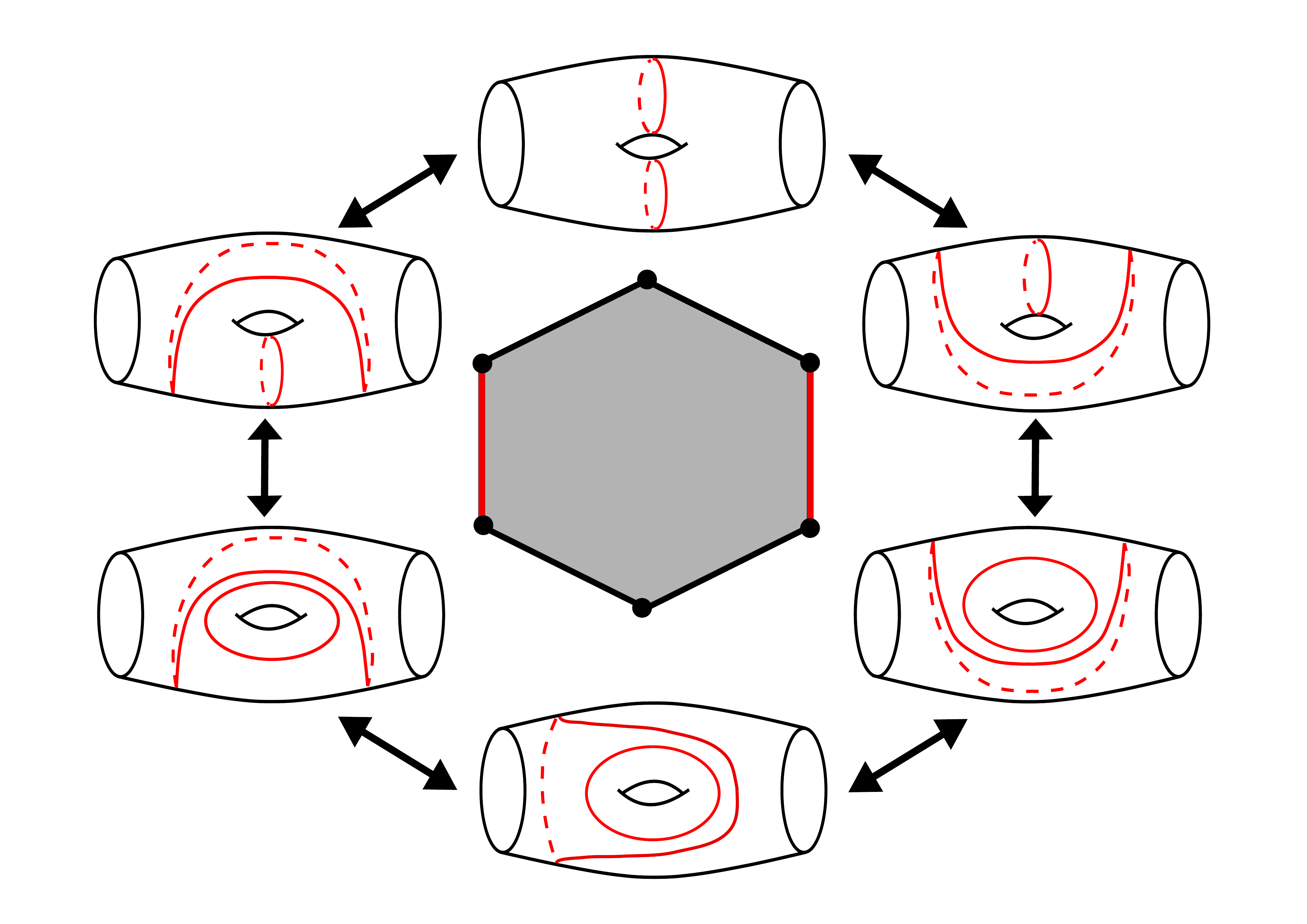}
    \caption{A 6AS cycle which is filled in by a hexagon. A-moves are labeled by black curves while S-moves are labeled by red curves.}
    \label{fig:6AS}
\end{figure}

\begin{figure}
    \centering
    \labellist
		\pinlabel {\LARGE{\textbf{$4S$}}} at 605 320
	\endlabellist
        \includegraphics[scale=.33]{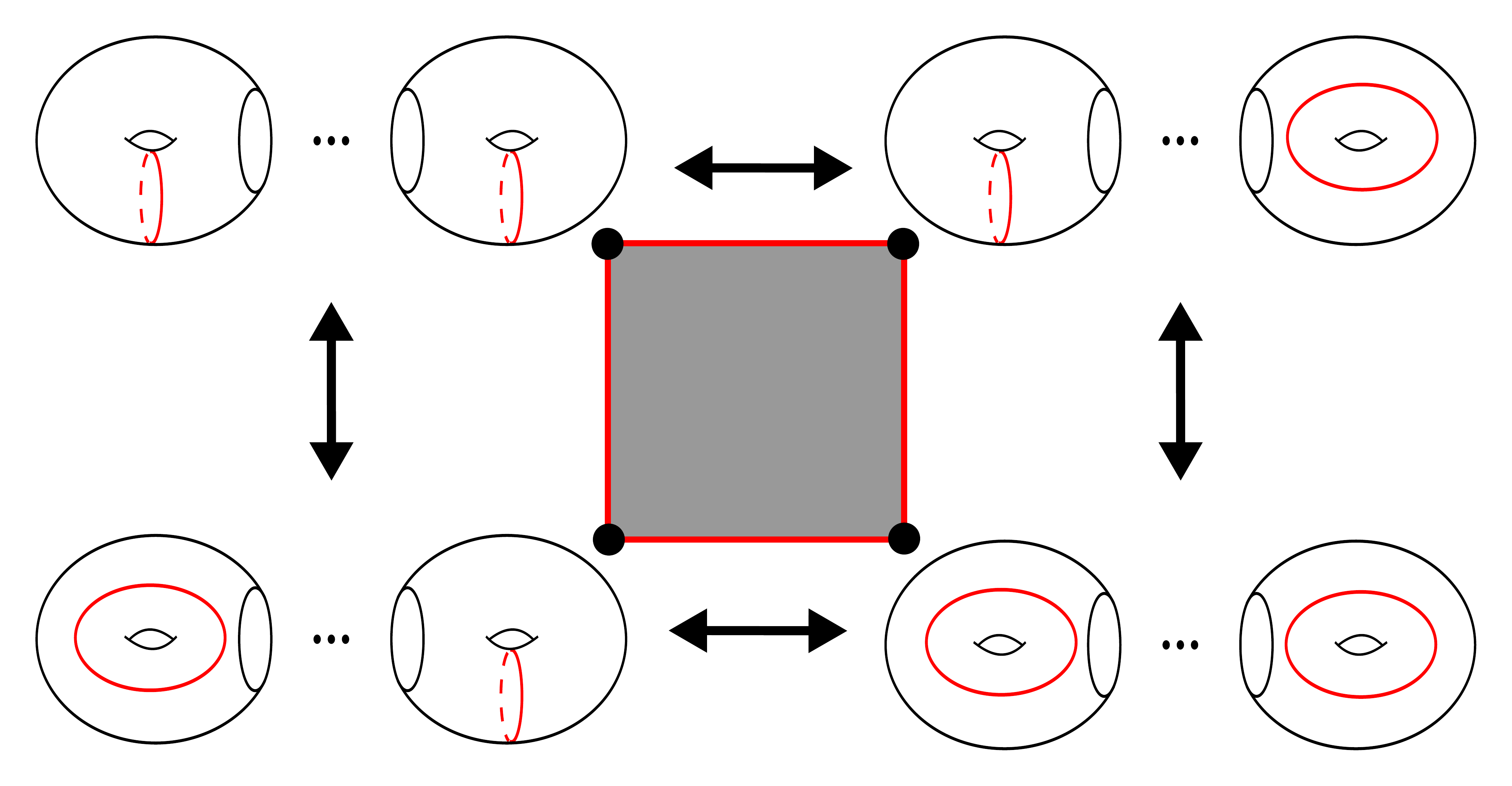}
    \caption{A 4S cycle of S-moves which is filled in by a square. There are also corresponding squares for disjointly supported A-moves and disjointly supported A- and S- moves.  We will denote these by 4S-squares, 4A-squares, and 4AS-squares, respectively.}
    \label{fig:4C}
\end{figure}

If $P$ is a pants decomposition for $\Sigma$, we may obtain a 3-dimensional handlebody with boundary $\Sigma$, by taking $\Sigma \times I$, attaching 2-handles to the curves in $P$, and then capping off any remaining spheres components with 3-handles. We denote this handlebody by $\mathcal{H}(P)$,   We say two handlebodies $H_1$ and $H_2$ with $\partial H_1 = \partial H_2 = \Sigma$ are equal if the identity map on the boundary extends to a homeomorphism from $H_1$ to $H_2$.  Note that there are inequivalent pants decompositions that produce equal handlebodies, namely, any two pants decompositions related by A-moves define the same handlebody.

Given a handlebody $H$ with $\partial H = \Sigma$, the $\textbf{handlebody set}$ of $H$ in $\mathcal{P}(\Sigma)$ is the set of pants decompositions $P$ of $\Sigma$ such that $\mathcal{H}(P) = H$. By our previous discussion, two pants decompositions related by an A-move lie in the same handlebody set. The following result originally proved by Luo shows that handlebody sets are in fact connected by $A$-moves.

\begin{lemma}(Corollary 1 of \cite{FeLu}) 
\label{connected}
Given a handlebody $H$ with $\partial H = \Sigma$ and two pants decompositions $P_1$ and $P_2$ in the handlebody set of $H$, there exists a path in $\mathcal{P}(\Sigma)$ consisting of exclusively A-edges between $P_1$ and $P_2$.  
\end{lemma}

 A \textbf{cut system} on a genus $g$ closed orientable surface $\Sigma$ is a set $\alpha = \{ \alpha_1,...,\alpha_g \}$ of pairwise disjoint non-separating essential curves on $\Sigma$.  At times it will be more natural to consider cut systems instead of pants decompositions. The following two lemmas will allow us to pass between these decompositions freely.

\begin{lemma}(Lemma 5 of \cite{JJ})
\label{cut system}
Any pants decomposition for $\Sigma$ contains at least $g$ non-separating curves.  Any choice of $g$ non-separating curves in a pants decomposition give a cut system of $\Sigma$.  
\end{lemma}

Given a pants decomposition, choosing a cut system via the previous lemma gives sufficient data to determine a handlebody. This is often enough information for the topological constructions in this paper. To return to the simplicial constructions of the pants complex we will need to complete cut systems or, more generally, any set of curves to a pants decomposition. The following lemma, which may proved by induction on the genus and an Euler characteristic argument, will allow us to do so.

\begin{lemma} 
\label{extending to pants decomposition}
	Let $H$ be a genus $g$ handlebody with boudnary $\Sigma$. Let  $k < 3g-3$, be a natural number and $\{c_1,...,c_k\}$ be a set of non-isotopic simple closed curves on $\Sigma$ such that $c_1,...,c_k$ all bound disjoint properly embedded disks in $H$.  Then there exist some additional simple closed curves on $\Sigma$, $\{c_{k+1},..., c_{3g-3} \}$, that bound disjoint properly embedded disks in $H$, so that $P = \{c_1,...,c_k,c_{k+1},...,c_{3g-3} \}$ is a pants decomposition for $\Sigma$ with $\mathcal{H}(P) = H$.  
\end{lemma}

\section{The 1-skeleton of the pants complex}

Given two vertices $P_1$ and $P_2$ in the pants complex, we obtain a Heegaard splitting of a closed orientable 3-manifold $\mathcal{M}^3(P_1, P_2) = \mathcal{H}(P_1) \cup_\Sigma \mathcal{H}(P_2)$.  Further, if we consider an ordering, as above with $P_1$ first then $P_2$, and if $\Sigma$ is oriented, then we can orient $\mathcal{H}(P_1)$ to agree with $\Sigma$ and $\mathcal{H}(P_2)$ to disagree with $\Sigma$, we then obtain an orientation of $\mathcal{M}^3(P_1, P_2)$, and we will henceforth assume that it carries this orientation.  Note that changing the orientation of $\Sigma$, or the order of $P_1$ or $P_2$, will change the orientation of the resulting 3-manifold.

Given an oriented edge $e$ between two vertices $P_1$ and $P_2$ in the pants complex, we may define a compact orientable 4-manifold $\mathcal{X}^4(W)$ with $\partial \mathcal{X}^4(e) = \mathcal{M}^3(P_1, P_2)$. We first need to determine the manifold $\mathcal{M}^3(P_1, P_2)$. By inspecting Heegaard diagrams of the manifolds involved we may determine the 3-manifolds associated to vertices connected by an edge. In particular, we obtain the following lemma.

\begin{lemma} 
	Let $P_1$ and $P_2$ be two pants decompositions of a surface $\Sigma$ of genus $g$.  If $P_1$ and $P_2$ are connected by an A-edge, then $\mathcal{H}(P_1) = \mathcal{H}(P_2)$ and therefore, $\mathcal{M}^3( P_1, P_2 ) = \sharp^g (S^1 \times S^2)$.  If $P_1$ and $P_2$ are connected by an S-edge, then $\mathcal{M}^3( P_1, P_2 ) = \sharp^{g-1} (S^1 \times S^2)$.  
\end{lemma}

Having determined the 3-manifold associated to an edge, we next define a unique oriented 4-manifold filling. If $\mathcal{X}^4(e)$ is an A-edge, then we fill the resulting $\sharp^g(S^1 \times S^2)$ with $\natural^g(S^1 \times D^3)$ . Similarly, if the edge is an S-edge, we fill the resulting $\sharp^{g-1}(S^1 \times S^2)$ with $\natural^{g-1}(S^1 \times D^3)$. By \cite{LP}, these fillings are unique.

Given an oriented walk, $W$, of arbitrary length, we construct $\mathcal{X}^4(W)$ by first constructing all of the 4-manifolds associated to each of the edges in $W$, and then gluing together each successive pair of edges along the common 3-dimensional handlebody via the identity map along their shared handlebody (see Figure \ref{fig:GluingTwoWedges}).  As a convention, if $W$ is just a single vertex $P$ then we just take the filling of $\mathcal{H}(P) \cup_\Sigma \mathcal{H}(P)$ by $\natural^g (S^1 \times B^3)$.

\begin{figure}
    \centering
    \labellist
		\pinlabel {\LARGE{\textbf{$v_1$}}} at 50 510
		\pinlabel {\LARGE{\textbf{$v_2$}}} at 420 510
		\pinlabel {\Large{\textbf{\textcolor{red}{$H_{v_1}$}}}} at 100 70
		\pinlabel {\Large{\textbf{\textcolor{blue}{$H_{v_2}$}}}} at 370 70
		\pinlabel{{$\natural S^1 \times B^3$}} at 235 160
	\endlabellist
        \includegraphics[scale=.22]{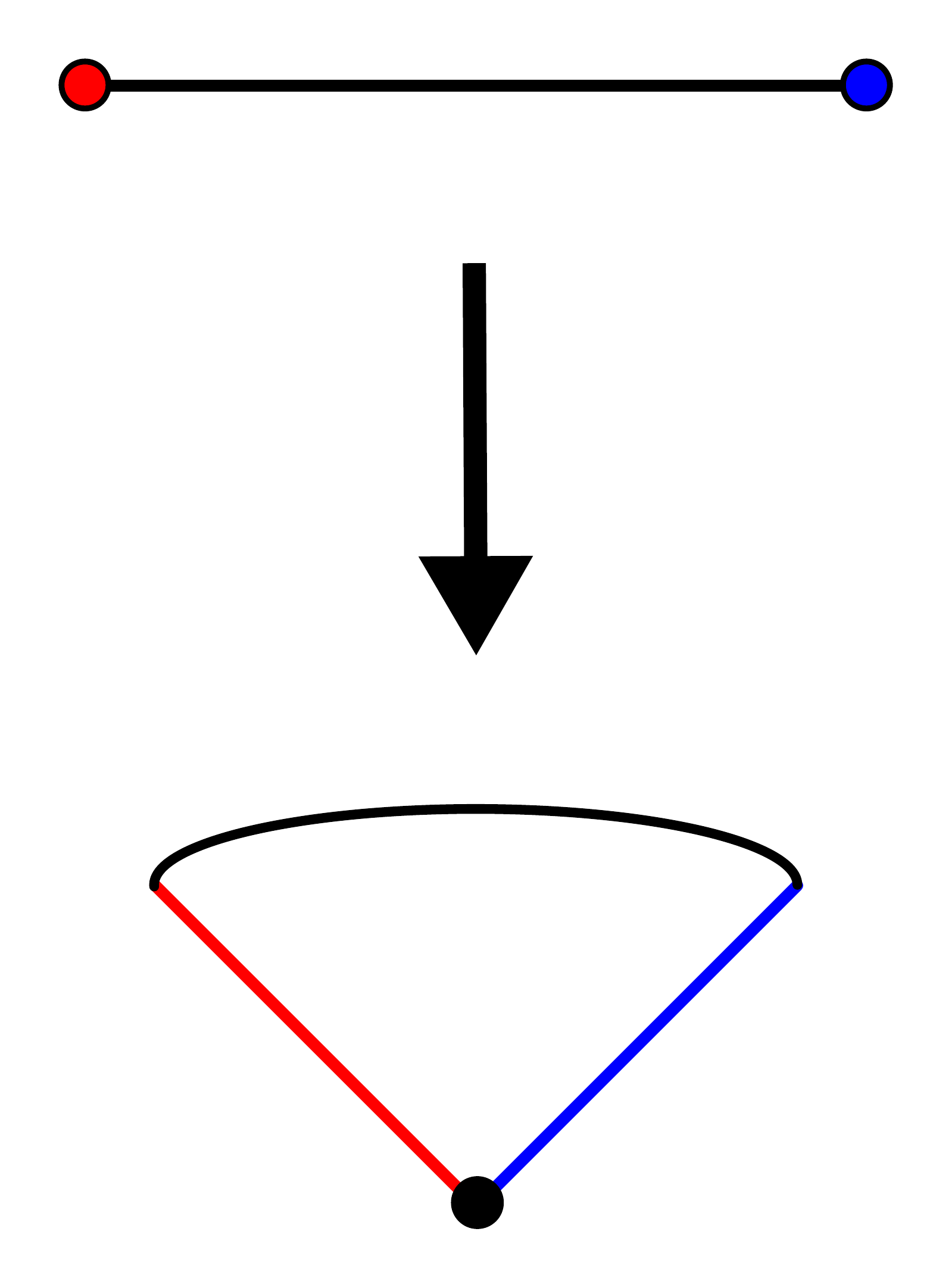}
    \caption{An edge in $P(\Sigma)$ corresponds to a connected sum of copies $S^1 \times S^2$ which may be filled in uniquely with the appropriate boundary sum of copies of $S^1 \times D^3$.}
    \label{fig:FillingInWedge}
\end{figure}
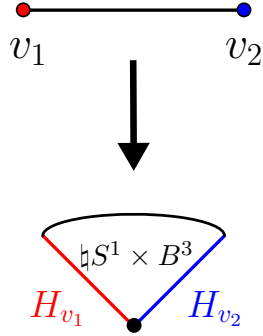

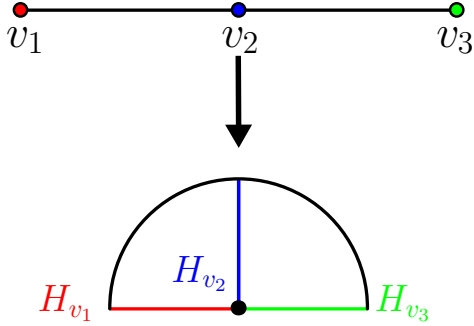
\begin{figure}
    \centering
    \labellist
		\pinlabel {\LARGE{\textbf{$v_1$}}} at 50 500
		\pinlabel {\LARGE{\textbf{$v_2$}}} at 420 500
		\pinlabel {\LARGE{\textbf{$v_3$}}} at 790 500
		
		\pinlabel {\Large{\textbf{\textcolor{red}{$H_{v_1}$}}}} at 120 50
		\pinlabel {\Large{\textbf{\textcolor{blue}{$H_{v_2}$}}}} at 350 100
		\pinlabel {\Large{\textbf{\textcolor{green}{$H_{v_3}$}}}} at 700 50
		
		\endlabellist
        \includegraphics[scale=.22]{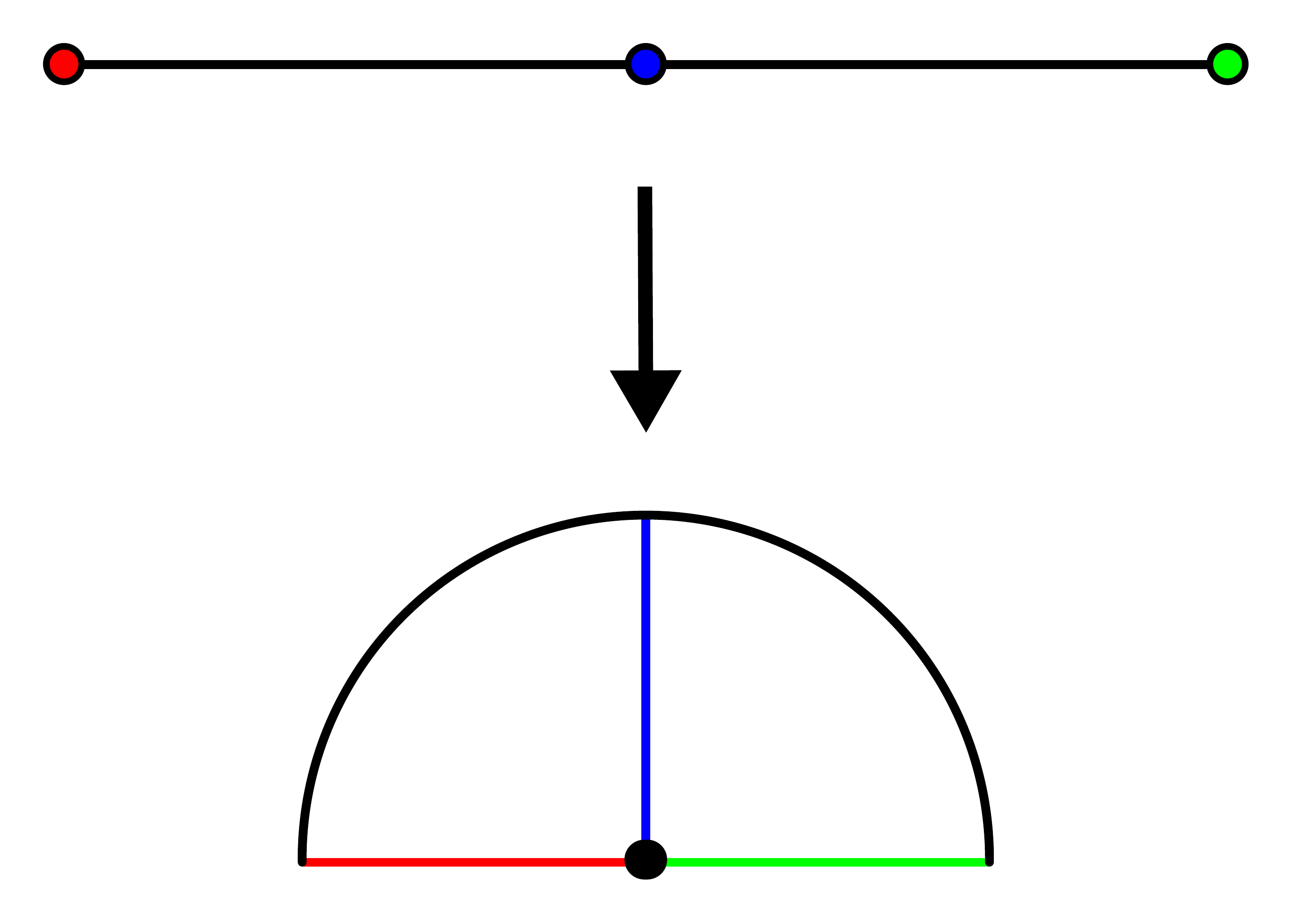}
    \caption{A path in $P(\Sigma)$ gives rise to the manifold obtained by gluing together the wedges corresponding to edges along the handlebodies of the shared vertices}
    \label{fig:GluingTwoWedges}
\end{figure}

We now address orientations. Assume from now on that we have a fixed orientation on $\Sigma$ and our walk $W$ has a fixed orientation, we will obtain an orientation on $\mathcal{X}^4(W)$ as follows:  for a single oriented edge from a vertex $P_1$ to a vertex $P_2$, orient $\mathcal{H}(P_1)$ with the orientation that induces the orientation on $\Sigma$ and orient $\mathcal{H}(P_2)$ with the orientation that induces the opposite orientation on $\Sigma$.  Then we can glue $\mathcal{H}(P_1)$ to $-\mathcal{H}(P_2)$ via the identity map and obtain an orientation on $\mathcal{M}^3(P_1 , P_2)$ that then induces an orientation on $\mathcal{X}^4(W)$.  If $W$ consists of multiple edges, then we can orient $\mathcal{X}^4(W)$ by orienting the wedges as above, and since each non-end vertex of $W$ has an edge coming in and an edge going out, the resulting identity maps between the handlebodies will be orientation-reversing and therefore we will obtain an orientation on all of $\mathcal{X}^4(W)$.  Note that, either switching the orientation of $\Sigma$ or switching the direction of $W$ will change the orientation on $\mathcal{X}^4(W)$.

We next seek to obtain a handle description for 4-manifolds given by a walk in the pants complex. We first need the following lemma. 

\begin{lemma} \label{lem:HandlebodyDualSurgery}
    Let $H$ be a handlebody and let $\gamma \subset \partial H$ be a curve  such that, for some properly embedded disk $D \subset H$, $|\gamma \cap D| = 1$. Then the result of pushing $\gamma$ into $H$, and doing surgery on $\gamma$ is again a handlebody. Moreover, if we do surgery on $\gamma$ using the surface framing, then  $\gamma$ bounds a disk in the surgered handlebody.
\end{lemma}

\begin{proof}

One way to prove that a 3-manifold is a handlebody is to find a collection of disjoint properly embedded disks which cuts the manifold into balls. Let $D_1,...,D_g$ be a collection of such disks for $H$ with $D_1 = D$. By sliding all of the other disks over $D_1$, we may arrange so that none of these disks intersect $\gamma$. We can then still cut along $D_2,...,D_{g}$ in the surgered manifold. We then only need to analyze what is happening in the solid torus $H \backslash \{D_2,...D_g \}$ containing the curve $\gamma$ where the surgery is occurring along with its dual disk $D_1$. 

Since $|\gamma \cap D_1| = 1$, and $\gamma$ is isotopic into the boundary, $\gamma$ is isotopic to the core curve $S^1 \times \{0\} \subset S^1 \times D^2$.  But any surgery on the core curve in a solid torus results again in a solid torus - one way to see this is that the part of the solid torus that is not affected by the surgery is just a collar neighborhood of the boundary.  If we give the surgery curve the surface framing, then there is a disk $D_1'$, disjoint from $D_2,...,D_{g}$, with $\partial D_1' = \gamma$, formed by taking the surgery disk for the push-in of $\gamma$ and extending it to the boundary by adding the annulus coming from the push-in process. We therefore conclude that $D_1', D_2, ..., D_g$ is a collection of properly embedded disks which cut the surgered manifold into a ball, and hence it is a handlebody.
\end{proof}

The following will used repeatedly to identify the 4-manifolds corresponding to paths and loops in the pants complex.

\begin{lemma} \label{handlebody decomposition}
	Let $W$ be an oriented walk in $\mathcal{P}(\Sigma)$ starting at $P_1$. The following process produces a handle decomposition of $\mathcal{X}^4(W).$ Start with $\mathcal{H}(P_1) \times I$; these are the 0- and 1-handles. For every (directed) S-edge in $W$, we take the new curve in the latter vertex, push it into  $\mathcal{H}(P_1) \times \{0\} \subset \mathcal{H}(P_1) \times I$, and give this curve the surface framing from $\Sigma =  \partial \mathcal{H}(P_1) \times \{0\}$. The curves that are seen later in the walk along $W$ are not pushed as far inside of $\mathcal{H}(P_1) \times \{0\}$ as earlier curves. These framed curves are the attaching curves for the 2-handles. There are no 3-handles or 4-handles.
\end{lemma}

\begin{proof}

    We start in the case where $W$ is just a single edge.  In the case where $W$ is an A-edge, then $\mathcal{X}^4(W) = \natural^g (S^1 \times D^3)$ and indeed this is the manifold that we obtain from our handlebody description, since no 2-handles are added.  
    
    In the case where $W$ is an S-edge, then $\mathcal{X}^4(W) = \natural^{g-1}(S^1 \times D^3)$ and so we must verify that the attaching sphere of the 2-handle that we are adding intersects the belt sphere of one of the 1-handles in exactly one point, and therefore the handles cancel giving the desired result.  Using Lemma \ref{cut system}, we can choose $g$ nonseparating curves in $P_1$ on $\Sigma$ that  form a cut system for $\mathcal{H}(P_1)$.  The belt spheres of the 1-handles are exactly  these $g$ nonseparating curves in $\Sigma$ together with the disks  on both sides of $\gamma \subset \sharp^g (S^1 \times S^2) = \partial (S^1 \times D^3)$, and by the definition of an S-move, and the convention for attaching a 2-handle stated in the lemma, we see that the attaching circle for the 2-handle intersects the belt sphere in exactly one point.  

	Finally, consider the case of a general walk $W = (P_1, P_2,... P_{n+1})$.  Since A-moves do not affect the resulting 4-manifold, we proceed by induction on the  number of S-moves that $W$ contains and we assume that ($P_n, P_{n+1}$) is an S-edge.  The base case was just discussed.  Let $W' = (P_1, ..., P_n)$.   By the inductive hypothesis, the 4-manifold $\mathcal{X}^4(W')$ has a handlebody diagram as described in the statement of the lemma.  Let $H = \mathcal{H}(P_{n})$.  Attach one end of $H \times I$ to this $H$ in the boundary to obtain a space that (after rounding corners) is still just $\mathcal{X}^4(W')$. Now attach a 2-handle along the new curve in $P_{n+1}$ framed by $\Sigma$ to the free end of $H \times I$ as in the statement of the lemma.

	We now verify that the union of $H \times I$ and this new 2-handle is indeed $\natural^{g-1} (S^1 \times D^3)$.  Using Lemma \ref{cut system}, we obtain a cut system for $H$ containing the curve that is changed by the S-move, such that the curves in the cut system bound disjoint disks in $H$.  These disks considered in both ends of $H \times I$ together with the curves cross $I$ form a set of belt spheres for the genus $g$ 4-dimensional handlebody $H \times I$ and the attaching circle for the 2-handle intersects exactly the disk bounding the curve corresponding to the S-move, and in exactly one point, thus verifying that we have $\natural^{g-1} (S^1 \times D^3)$ as desired.  Furthermore, by Lemma \ref{lem:HandlebodyDualSurgery}, when we look at the two handlebodies in the boundary of $H \times I$ together with this 2-handle, we have exactly $\mathcal{H}(P_n)$ and $\mathcal{H}(P_{n+1})$. Therefore, what we have attached is a 4-dimensional filling of the desired handlebodies by $\natural^{g-1} (S^1 \times D^3)$, which, by $\cite{LP}$, can only be done in one way. 
\end{proof}

\begin{figure}
    \centering

        \includegraphics[scale=.22]{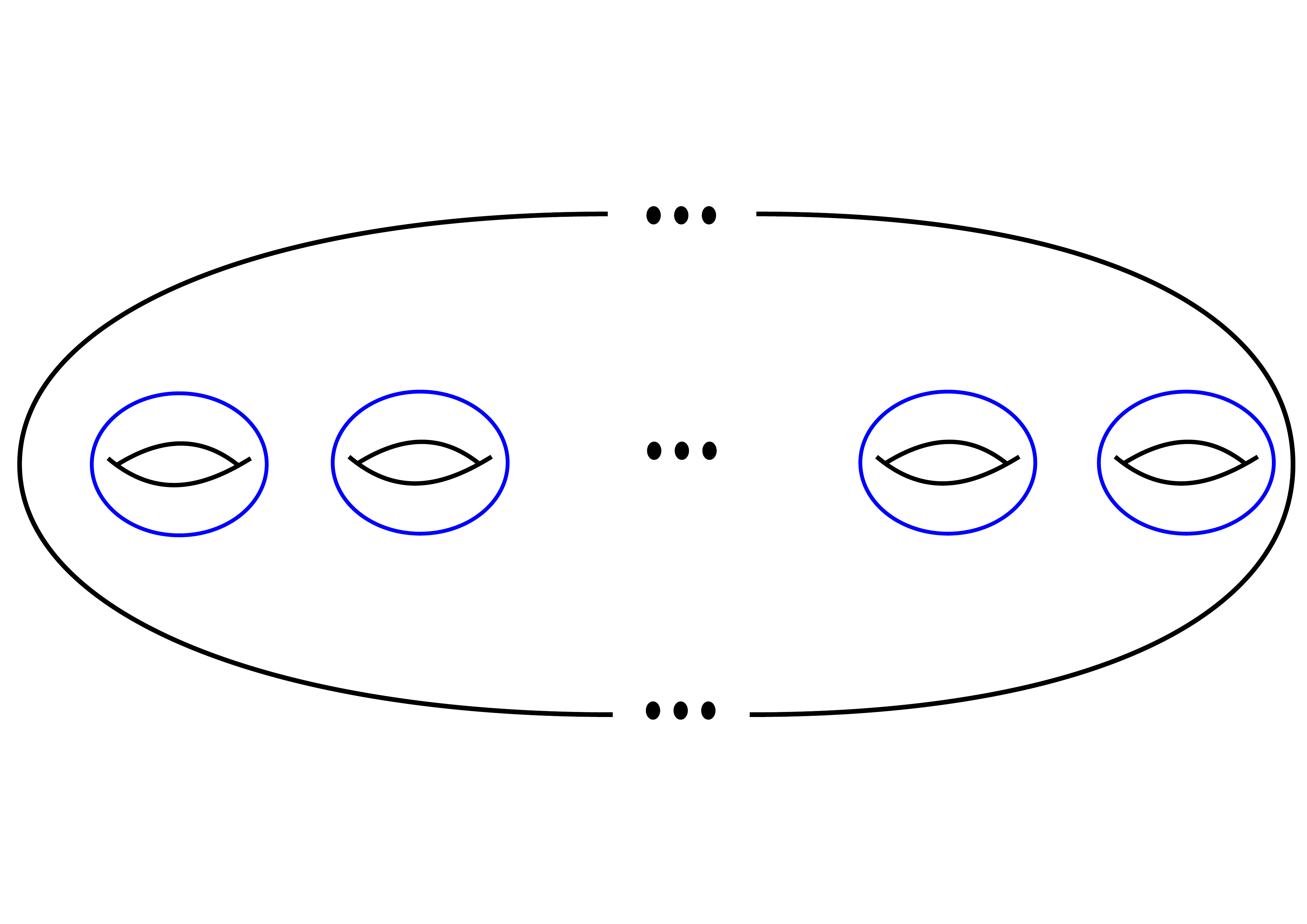}
    \caption{A cut system for the standard handlebody. These curves become the 1-handles in the dotted circle notation for 1-handles in a Kirby diagram. }
    \label{fig:betaHandlebody}
\end{figure}

 Since the mapping class group acts transitively on the set of handlebodies with a given boundary, we may apply a mapping class, and insert or delete A-moves, to assume that a walk starts in a pants decomposition which contains the cut system shown in Figure $\ref{fig:betaHandlebody}$. We may then use Lemma \ref{handlebody decomposition}, to obtain a Kirby diagram for the manifold $X^4(W)$. Namely, the cut system for the handlebody shown in Figure \ref{fig:betaHandlebody} become dotted circles, representing the 1-handles. S-moves give rise to 2-handles as in the previous lemma; these are pushed towards the ``outside" of the surface pictured.

The following result uses Waldhausen's theorem \cite{FW} on Heegaard splittings of connect sums of $S^1 \times S^2$, which is the primary 3-dimensional result that we will make use of.

\begin{lemma} \label{Waldhausen}
Let $P_1$ and $P_2$ be vertices in $\mathcal{P}(\Sigma_g)$, with $\mathcal{M}^3(P_1,P_2) \cong \sharp^k(S_1 \times S_2)$.  Then there exists a walk $W$ in $\mathcal{P}(\Sigma)$ with $\mathcal{X}^4(W) \cong \natural^k(S^1 \times D^3)$.  
\end{lemma}

\begin{proof}
Note that $P_1$ and $P_2$ form a genus $g$ Heegaard diagram for $\sharp^m (S^1 \times S^2)$. As an immediate consequence of Waldhausen's theorem \cite{FW}, there exist vertices $P_1'$ and $P_2'$ in the same handlebody sets as $P_1$ and $P_2$, respectively, which are standard, in that there are $k$ parallel sets of non-separating curves and $g-k$ dual sets of non-separating curves in these pants decompositons. The path of length $g-k$ between these pants decompositions which turns each curve in $P_1$ to its corresponding dual in $P_2$ gives rise to a handle decomposition consisting of $g$ 1-handles and $g-k$ 2-handles. These 2-handles are dual to 1-handles, so the resulting manifold has a handle decomposition consisting of just $k$ 1-handles, and so is diffeomorphic to $\natural^k(S^1 \times D^3)$.
\end{proof}

We can also construct a closed 4-manifold given a loop $L$ in $\mathcal{P}(\Sigma)$ just as in the above construction, but when the loop returns to the vertex that we start on, we glue the identical handlebodies together using the identity map.  Further, if the loop is oriented, the resulting 4-manifold obtains an orientation. We denote the resulting closed orientable manifold by $\mathcal{X}_C^4(L)$.  Note that changing the orientation of $L$ or $\Sigma$ changes the orientation of $\mathcal{X}_C^4(L)$.

We next seek to prove that every 4-manifold arises in this fashion. To do this, we will convert a handle decomposition into a loop in the pants complex - this is very similar to the handlebody proof that every 4-manifold admits a trisection \cite{GK}. Throughout the proof, the reader may find it helpful to consult Figures \ref{fig:Genus2CP2} and \ref{fig:4Cgenus2}.  These figures, read left to right, show how to obtain a handlebody diagram from a loop in the pants complex.  The following proof illustrates the reverse of this, so these figures, read right to left, provide examples of the procedure.  

\existence

\begin{figure}
    \centering
    \includegraphics[scale=.2]{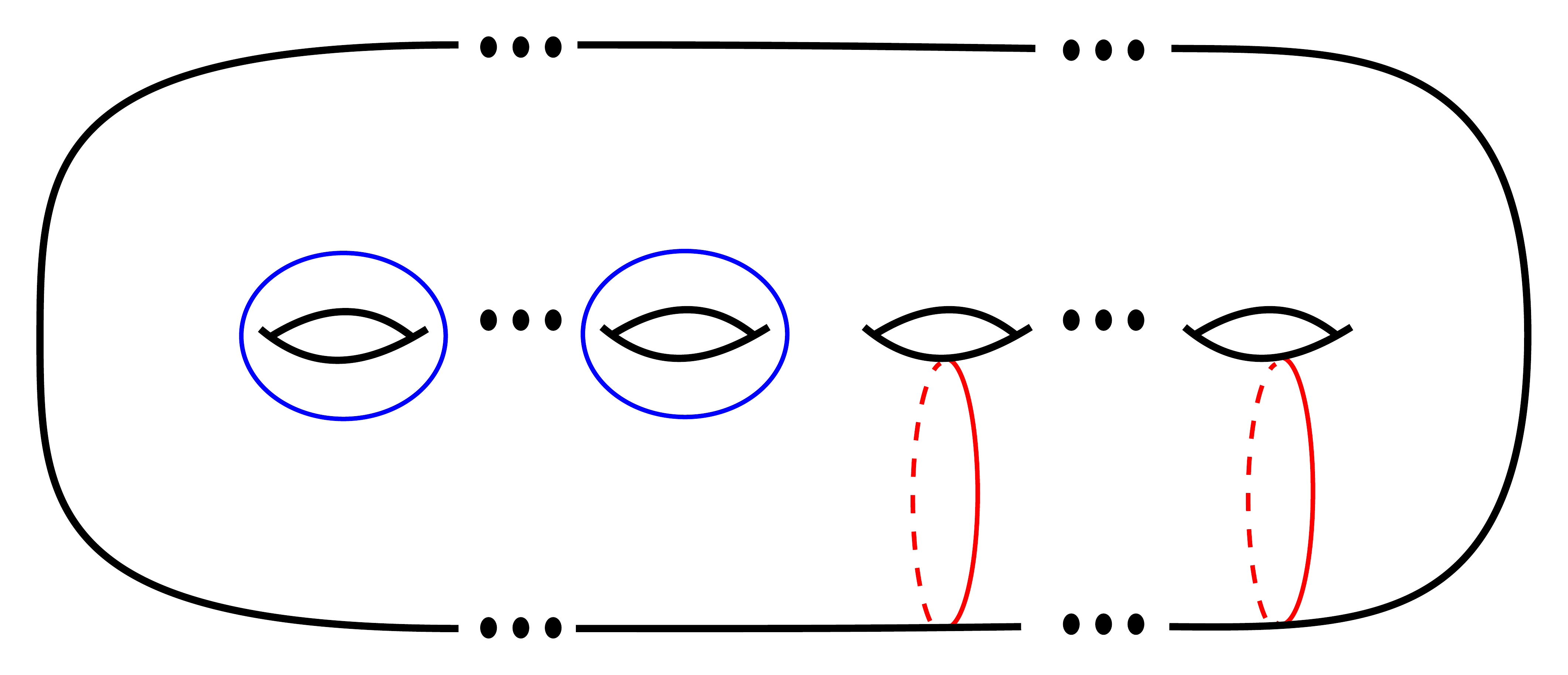}
    \caption{The handlebody $Q$ used in Thereom \ref{thm:existence}}
    \label{fig:QHandlebody}
\end{figure}

\begin{figure}
    \centering
    \includegraphics[scale=.5]{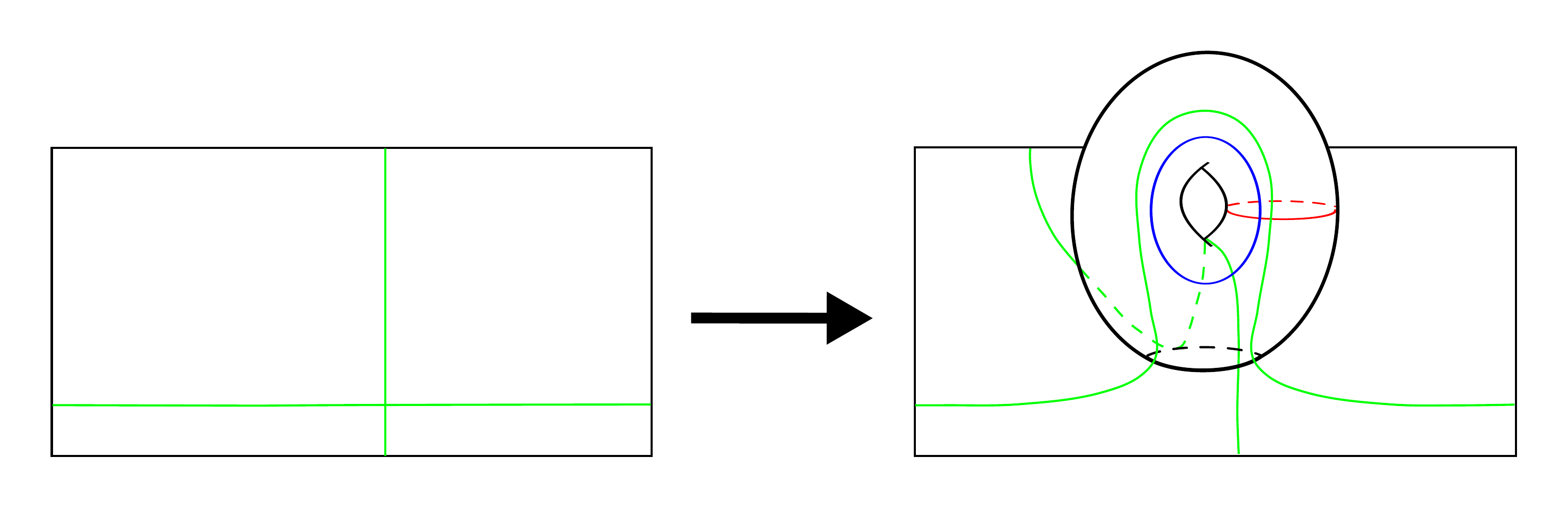}
    \caption{Stabilizing the Heegaard surface of the boundary of the 0- and 1-handles allows us to eliminate any intersections which arise between the 2-handles when they are projected onto the Heegaard surface.}
    \label{fig:ResolvingIntersection}
\end{figure}

\begin{figure}
    \centering
    \includegraphics[scale=.5]{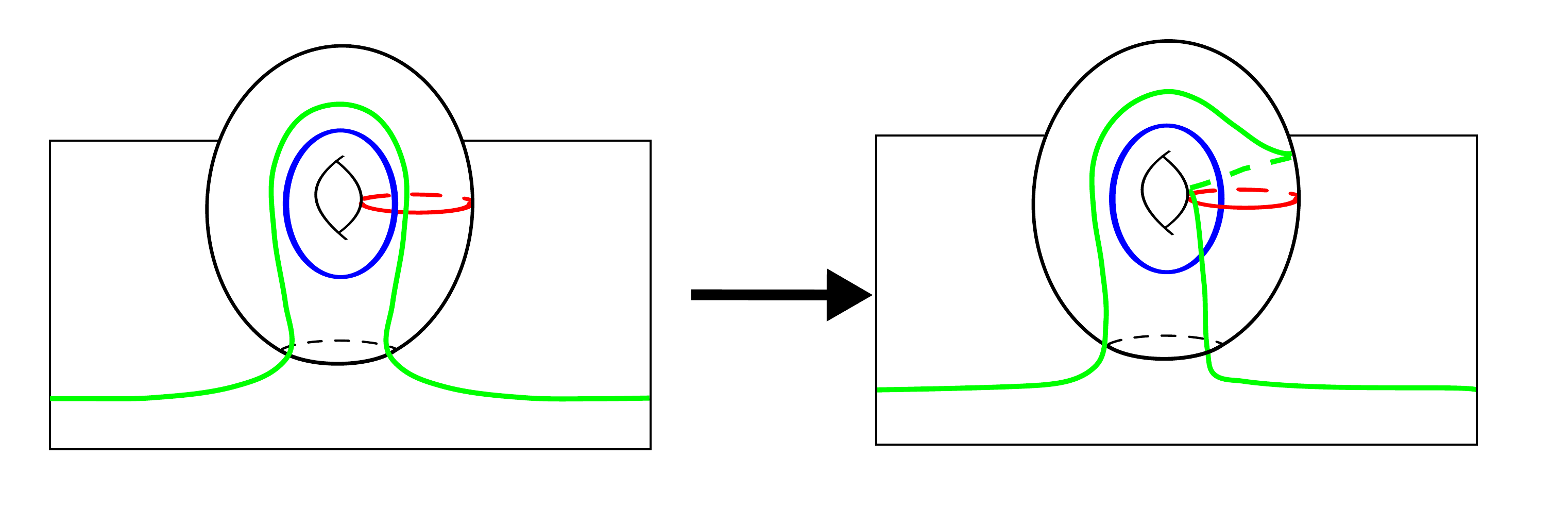}
    \caption{By twisting an attaching circle of a 2-handles around the corresponding dual $\alpha$ curve, we can make the handle framing match the surface framing.}
    \label{fig:CorrectingFraming}
\end{figure}

\begin{proof}

	Let $X^4$ be an arbitrary closed 4-manifold. Fix a handlebody decomposition diagram for $X$.  Let $\Sigma_0$ be a Heegaard splitting surface of the boundary of the 4-dimensional handlebody that we see after just attaching the 0- and 1-handles in the construction of $X$.  Let $l = l_1 \cup \cdots \cup l_n$ be the framed link that describes how the 2-handles are attached.  Project $l$ onto $\Sigma_0$ and then stabilize the Heegaard surface $\Sigma_0$ to obtain a new Heegaard surface $\Sigma$ in the following way.  First, stabilize $\Sigma_0$ in order to make the $l_i$ embedded  as in Figure \ref{fig:ResolvingIntersection}. If needed, stabilize further so that for each curve, $l_i$, there is a curve, $\alpha_i$, embedded in the surface so that $\alpha_i$ intersects $l$ in exactly one point.  Call this resulting surface $\Sigma$.  By twisting the $l_i$ around the $\alpha_i$ as in Figure \ref{fig:CorrectingFraming}, we can ensure that the framing on each $l_i$ is the same as the framing coming from the surface embedding, and we will assume that $l$ is sitting in $\Sigma$ in this way.  

	We now construct our loop $L$ in $\mathcal{P}(\Sigma)$ with $\mathcal{X}^4_C(L) \cong X$.  We will construct $L$ so that the handlebody decomposition of $\mathcal{X}_C^4(L)$ that we see from Lemma \ref{handlebody decomposition} is identifiable with the given handlebody decomposition of $X$.

	Suppose that $g$ is the genus of $\Sigma$ and $k$ is the genus of $\Sigma_0$ (i.e. the number of 1-handles in the given handlebody decomposition of $X$).  Our construction of $W$ will take place in a few stages. We start by constructing the 1-handles of $X^4$.  Take a pants decomposition of $\Sigma$ that contains the cut system in Figure \ref{fig:betaHandlebody}. By performing $g-k$ S-moves, we arrive at a pants decomposition that contains the cut system in Figure \ref{fig:QHandlebody}, which we call $Q$.  We call this walk $W_1$.  At this point, via Lemma \ref{handlebody decomposition}, we see that we have constructed a genus $k$ 4-dimensional 1-handlebody, and all of the 1-handles that have been cancelled are exactly the handles that do not appear in the given handlebody decomposition of $X$.  

	Let $n_i$ denote the boundary of the punctured torus that is a regular neighborhood of $\alpha_i \cup l_i$.  Note that the $n_i$ together with all of the $\alpha_i$ form a collection of disjoint simple closed curves that bound disjoint properly embedded disks in  $\mathcal{H}(Q)$.  Let $R$ be a pants decomposition obtained from applying Lemma \ref{extending to pants decomposition} to the union of the $n_i$ and $\alpha_j$, so that $\mathcal{H}(R) = \mathcal{H}(Q)$.  By Lemma \ref{connected}, we can get from $Q$ to $R$, using a walk, which we call $W_2$, with just A-moves.  Therefore, $\mathcal{X}^4(W_1W_2) \cong \mathcal{X}^4(W_1) \cong \natural^k(S^1 \times D^3)$, as no new handles have been added.  

	Now we are in position to attach the desired 2-handles by doing S-moves.  Namely, for each $\alpha_i$ curve do the S-move that turns $\alpha_i$ into $l_i$.  Let $W_3$ be the walk starting at the vertex $R$ that consists of this sequence of S-moves.  By Lemma \ref{handlebody decomposition}, we see that $\mathcal{X}^4(W_1W_2W_3)$ is diffeomorphic to the 0-, 1- and 2-handles in the handlebody decomposition of $X$. Since $X$ is a closed 4-manifold, we must have that the boundary of  $\mathcal{X}^4(W_1W_2W_3)$ is $\natural^m (S^1 \times D^3)$ for some $m$. Moreover, the first and last handlebodies of $W_1W_2W_3$  form a Heegaard splitting for the boundary. By Lemma \ref{Waldhausen},  there exists a walk $W_4$ from the end of $W_3$ to the beginning of $W_1$ with $\mathcal{X}^4(W_4) = \natural^m (S^1 \times D^3)$, and therefore, again since 3- and 4-handles glue in uniquely by \cite{LP}, we see that $\mathcal{X}^4_C(W_1W_2W_3W_4) \cong X$.

\end{proof}

\section{The 2-skeleton of the pants complex and cobordisms of 4-manifolds}
\label{sec:2-cells}
A homotopy in a 2-complex can be seen as replacing parts of the boundary of a 2-cell by the rest of the boundary, as in Figure \ref{fig:HomotopingPath}. We seek to understand how our 4-manifolds change under this operation. To this end, we must understand the 4-manifolds that arise along all connected subsets of the boundaries of the 2-cells.

Note that the mapping class group of $\Sigma$ acts on the set of handlebodies with boundary $\Sigma$, and that this action is transitive.  If $H$ is a handlebody with $\partial H = \Sigma$ and $\phi$ is a mapping class of $\Sigma$, we will denote the result of the action of $\phi$ on $H$ by $\phi \cdot H$.  

\begin{lemma} \label{mapping class}
	Let $W = P_1 \cdots P_n$, $W_2 = Q_1 \cdots Q_n$ be two walks in $\mathcal{P}(\Sigma)$, and let $\phi$ be a mapping class of $\phi$ such that $\phi \cdot \mathcal{H}(P_i) = \mathcal{H}(Q_i)$ for all $1 \leq i \leq n$.  Then $\mathcal{X}^(W_1) \cong \mathcal{X}^4(W_2)$.  The analogous statement also holds for closed loops.  
\end{lemma} 

\begin{proof}
	This is an immediate application of \cite{LP} applied to each of the wedge pieces.  
\end{proof}

We now go through an extended example, analyzing the 4-manifolds obtained from the boundary of one of the two-cells in $\mathcal{P}(\Sigma)$.  Many of the arguments will be repeated for the other 2-cells and this will form the core of the results that follow.  Lemma \ref{handlebody decomposition} will be applied throughout this example and the examples that follow in order to obtain handlebody decompositions so that we can recognise the relevant manifolds.  

\subsection{3S-Triangles}
\label{subsec:3S}
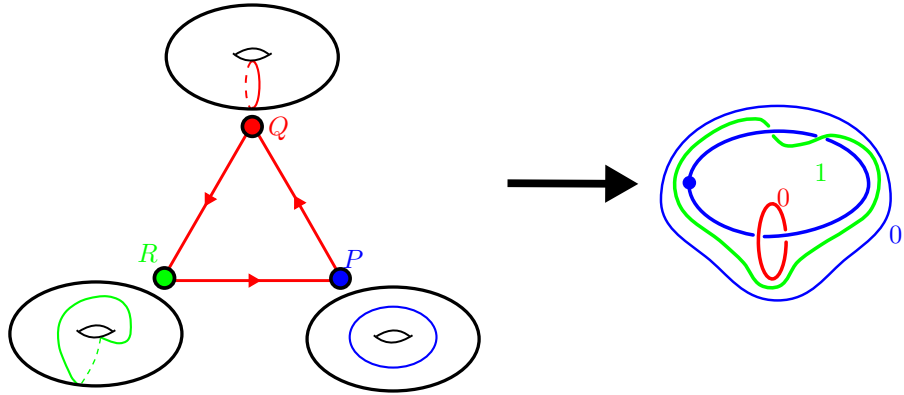
\begin{figure}
    \centering
    \labellist
		\pinlabel {\textcolor{red}{\small{\textbf{$0$}}}} at 650 190
		\pinlabel {\textcolor{blue}{\small{\textbf{$0$}}}} at 740 160
		\pinlabel  {\textcolor{green}{\small{\textbf{$1$}}}} at 680 210
		\pinlabel {\textcolor{blue}{\small{\textbf{$P$}}}} at 305 140
		\pinlabel {\textcolor{red}{\small{\textbf{$Q$}}}} at 245 245
		\pinlabel  {\textcolor{green}{\small{\textbf{$R$}}}} at 140 145
		
		\endlabellist
        \includegraphics[scale=.48]{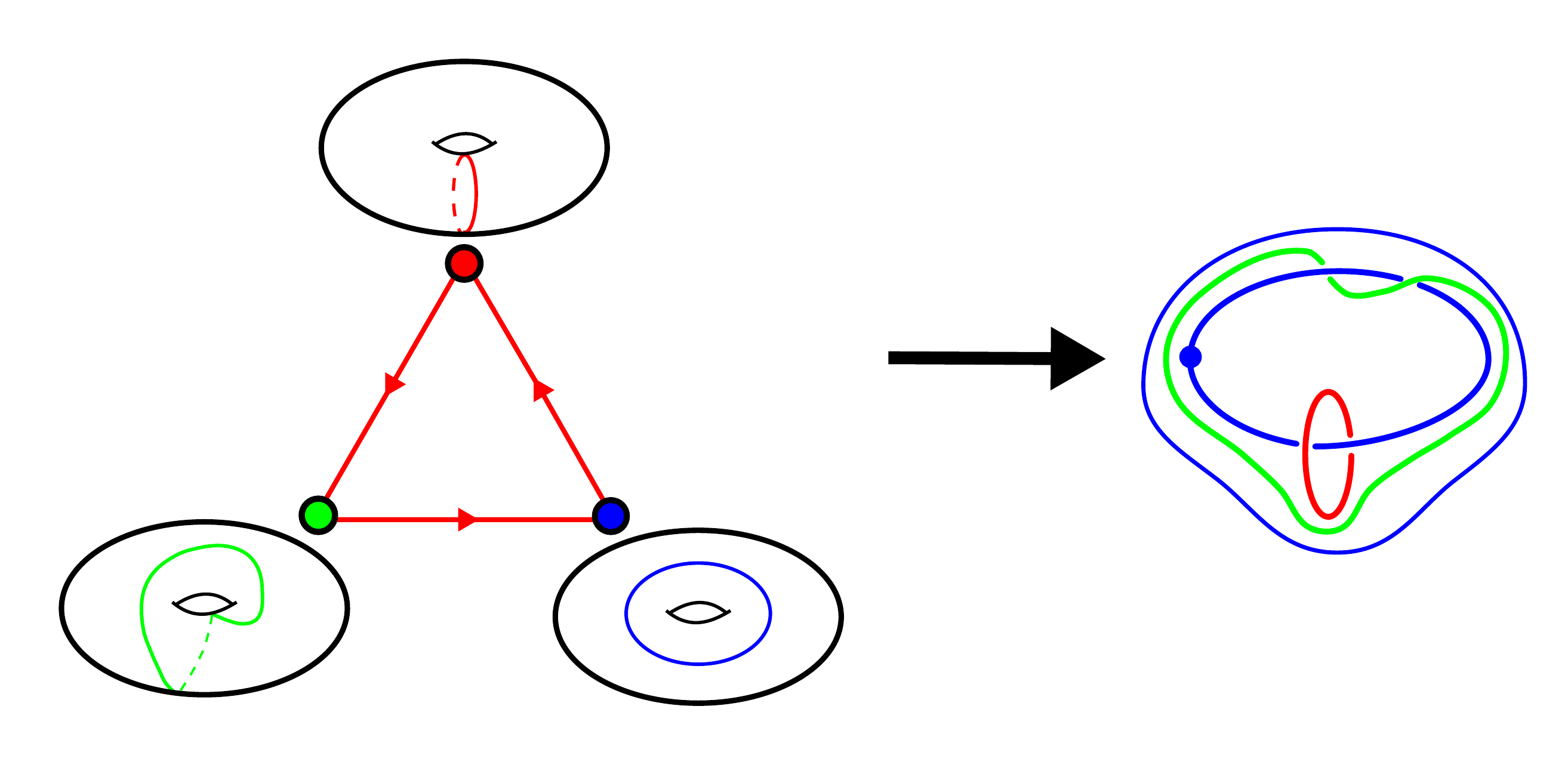}
    \caption{A 3S triangle in $P(\Sigma_1)$ gives rise to the Kirby diagram on the right. After sliding the green curve over the red curve and cancelling the red-blue pair, we see that this manifold is $\mathbb{C}P^2$.}
    \label{fig:Genus1CP2}
\end{figure}

Let $\Sigma$ be a torus and recall that the pants decomposition $\mathcal{P}(\Sigma)$ is defined to have vertices as isotopy classes of essential curves, only S-edges, and 3s-triangles as 2-cells.  Consider the oriented boundary of the 3S-triangle in Figure \ref{fig:Genus1CP2}. After a handle slide, it is evident that $\mathcal{X}^4_C(PQRP) = \mathbb{C}P^2$.  Note that if we instead traverse the triangle in the other direction, we obtain $\overline{\mathbb{C}P}^2$.  By applying an appropriate mapping class, we can interchange any of the pants decompositions in the triangle.  By Lemma \ref{mapping class},  all of the resulting segments are diffeomorphic to the corresponding segments in $PQRP$.  Thus, we have seen that any edge in the $PQRP$ triangles corresponds to a 4-ball, any pair of adjacent edges corresponds to $\mathbb{C}P^2 - \mathring{B}^4$ or $\overline{\mathbb{C}P}^2 - \mathring{B}^4$ depending on the orientation of the edges, and the whole triangle corresponds to $\mathbb{C}P^2$ or $\overline{\mathbb{C}P^2}$, again depending on the orientation.  Note that for any set of curves on $\Sigma$ that form a 3S-triangle $\Delta$, there is a mapping class $\phi$ of $\Sigma$ that sends the vertices of $\Delta$ to $P,Q,R$ in some order, and therefore the above analysis carries over for any such $\Delta$.  

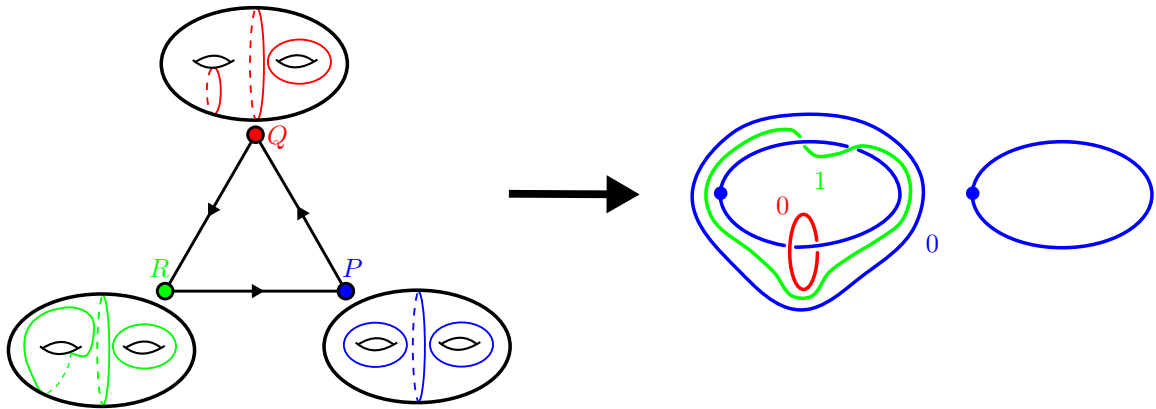
\begin{figure}
    \centering
    \labellist
		\pinlabel {\textcolor{red}{\small{\textbf{$0$}}}} at 650 190
		\pinlabel {\textcolor{blue}{\small{\textbf{$0$}}}} at 770 160
		\pinlabel  {\textcolor{green}{\small{\textbf{$1$}}}} at 680 210
		\pinlabel {\textcolor{blue}{\small{\textbf{$P$}}}} at 305 140
		\pinlabel {\textcolor{red}{\small{\textbf{$Q$}}}} at 245 245
		\pinlabel  {\textcolor{green}{\small{\textbf{$R$}}}} at 150 140
		
		\endlabellist
        \includegraphics[scale=.48]{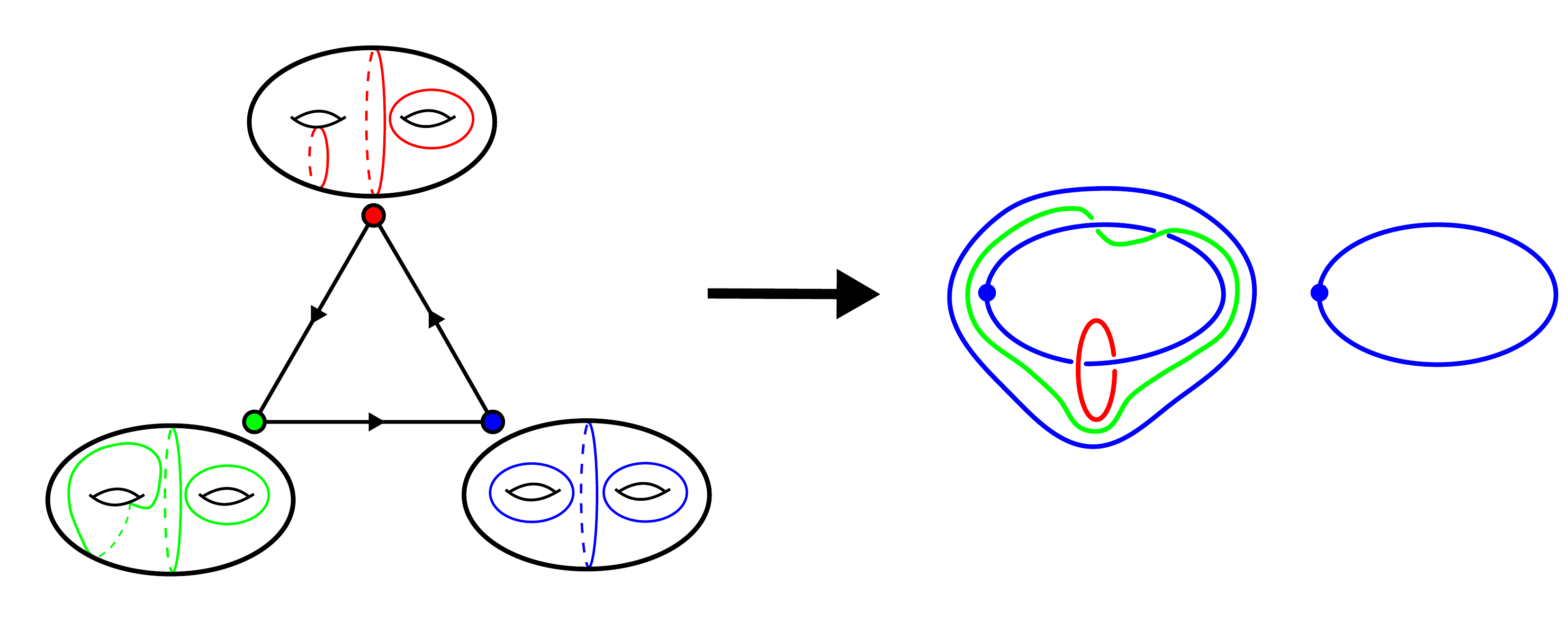}
    \caption{A 3S triangle in $P(\Sigma_2)$ gives rise to the Kirby diagram on the right. After sliding the green curve over the red curve and cancelling the red-blue pair, we see that this manifold is $\mathbb{C}P^2 \# (S^1 \times S^3)$.}
    \label{fig:Genus2CP2}
\end{figure}

An analogous analysis for 3S-triangles on higher genus surfaces is very similar.  For example, in a genus 2 surface $\Sigma$ in Figure \ref{fig:Genus2CP2} we have $\mathcal{X}^4_C(PQRP) = (S^1 \times D^3) \sharp \mathbb{C}P^2$.  As in the torus example, we see by Lemma \ref{mapping class} that this analysis holds for any individual edge in the triangle, and any pair of adjacent edges, and any 3S-triangle in $\mathcal{P}(\Sigma)$.  Similarly, this all goes through in the genus $g$ case where an edge will give $\natural^{g-1}(S^1 \times D^3)$, a pair of adjacent edges will give $\natural^{g-1}(S^1 \times D^3) \sharp \mathbb{C}P^2$ or $\natural^{g-1}(S^1 \times D^3) \sharp \overline{\mathbb{C}P}^2$ depending on the orientation, and the whole triangle will give $\sharp^g(S^1 \times S^3) \sharp \mathbb{C}P^2$ or $\sharp^g(S^1 \times S^3) \sharp \overline{\mathbb{C}P}^2$ depending on the orientation.  

\subsection{3A-Triangles, 4A-squares and 5A-Pentagons}
\label{subsec:A-Cells}
	The 3A-triangle and the 5A-pentagon both give rise to $\sharp^g (S^1 \times S^3)$ when $\Sigma$ has genus $g$, since, by Lemma \ref{handlebody decomposition}, the resulting manifold is built with a 0-handle, $g-$ 1-handles, $g$ 3-handles, and a 4-handles.  Further, again by Lemma \ref{handlebody decomposition} all of the edges and sequences of adjacent edges give rise to $\natural^g(S^1 \times D^3)$.  The same also holds true for the 4A-square.

	\begin{figure}
    \centering
   
    \labellist
    \pinlabel{P} at 322 120
    \pinlabel{Q} at 322 270
    \pinlabel{R} at 130 270
    \pinlabel{S} at 130 120
	\pinlabel {\textcolor{red}{\small{\textbf{$0$}}}} at 660 190
	\pinlabel {\textcolor{blue}{\small{\textbf{$0$}}}} at 630 265
	\pinlabel  {\textcolor{green}{\small{\textbf{$0$}}}} at 860 190
	\pinlabel  {\textcolor{purple}{\small{\textbf{$0$}}}} at 950 210
    \endlabellist
     \includegraphics[scale=.35]{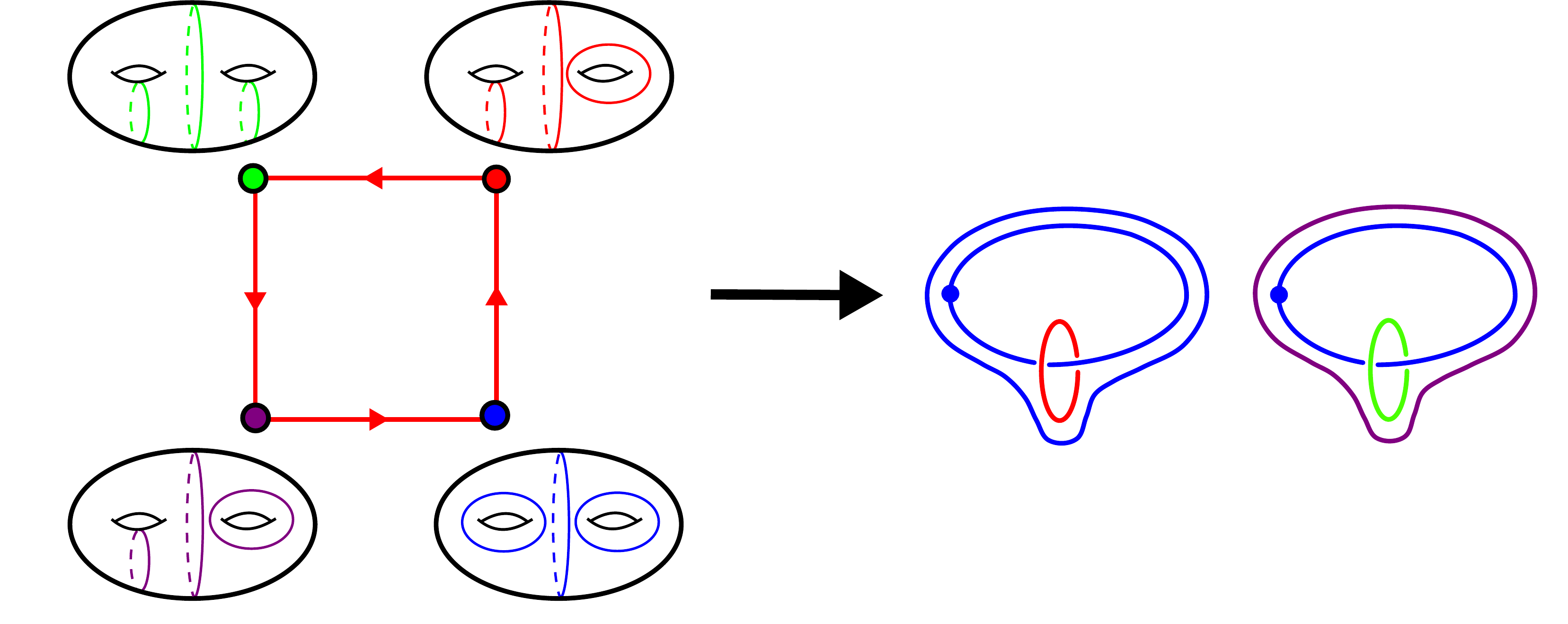}
    \caption{The 4S-square analyzed in Example \ref{subsec:4S}. In $P(\Sigma_2)$, this square gives rise to the Kirby diagram on the right. After cancelling the 1-2 pairs one recognizes this as a Kirby diagram for $S^4$ consisting of cancelling 2- and 3-handles.} 
    \label{fig:4Cgenus2}
\end{figure}
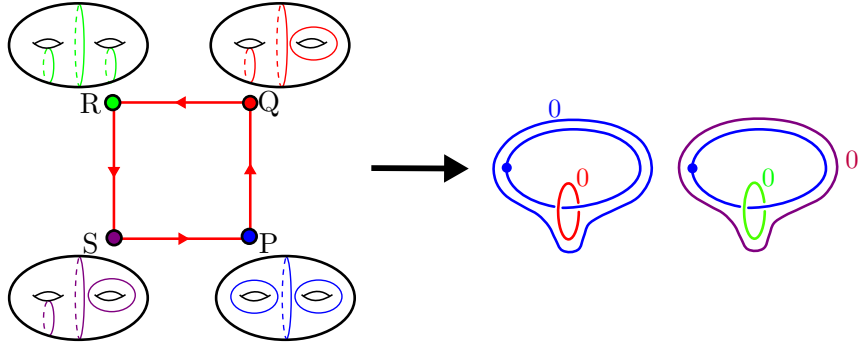

\subsection{4S-Squares}
\label{subsec:4S}
	In this example, we look at the 4S-squares.  We start with the genus 2 case shown in Figure \ref{fig:4Cgenus2}.  From the Kirby diagram, we see that $\mathcal{X}^4_C(PQRSP) = S^4$.  By the symmetry of the above diagrams, and Lemma \ref{mapping class}, we find that the above analysis holds for all of the edges of of the square, and all pairs and 3-tuples of adjacent edges. Again by Lemma \ref{mapping class}, this holds for all such 4S-squares in $\mathcal{P}(\Sigma)$.

	In the case of a genus $g$ surface $\Sigma$, as in the previous examples, when $g > 2$ this simply adds more 1-handles to the above situation.  So in this case the edges will result in $\natural^{g-1}(S^1 \times D^3)$, a pair of adjacent edges will give $\natural^{g-2}(S^1 \times D^3)$, three adjacent edges will give $\natural^{g-2}(S^1 \times D^3) \natural (S^2 \times D^2)$, and the whole square as a closed manifold will be $\sharp^{g-2}(S^1 \times S^3) \sharp (S^2 \times S^2)$.

\subsection{4AS-Squares and 6AS-Hexagons}
\label{subsec:4AS6AS}
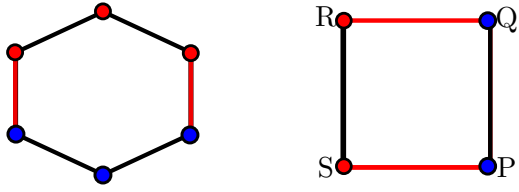
\begin{figure}
    \centering
   
    \labellist
    \pinlabel{P} at 657 50
    \pinlabel{Q} at 657 235
    \pinlabel{R} at 430 240
    \pinlabel{S} at 430 50
    \endlabellist
     \includegraphics[scale=.3]{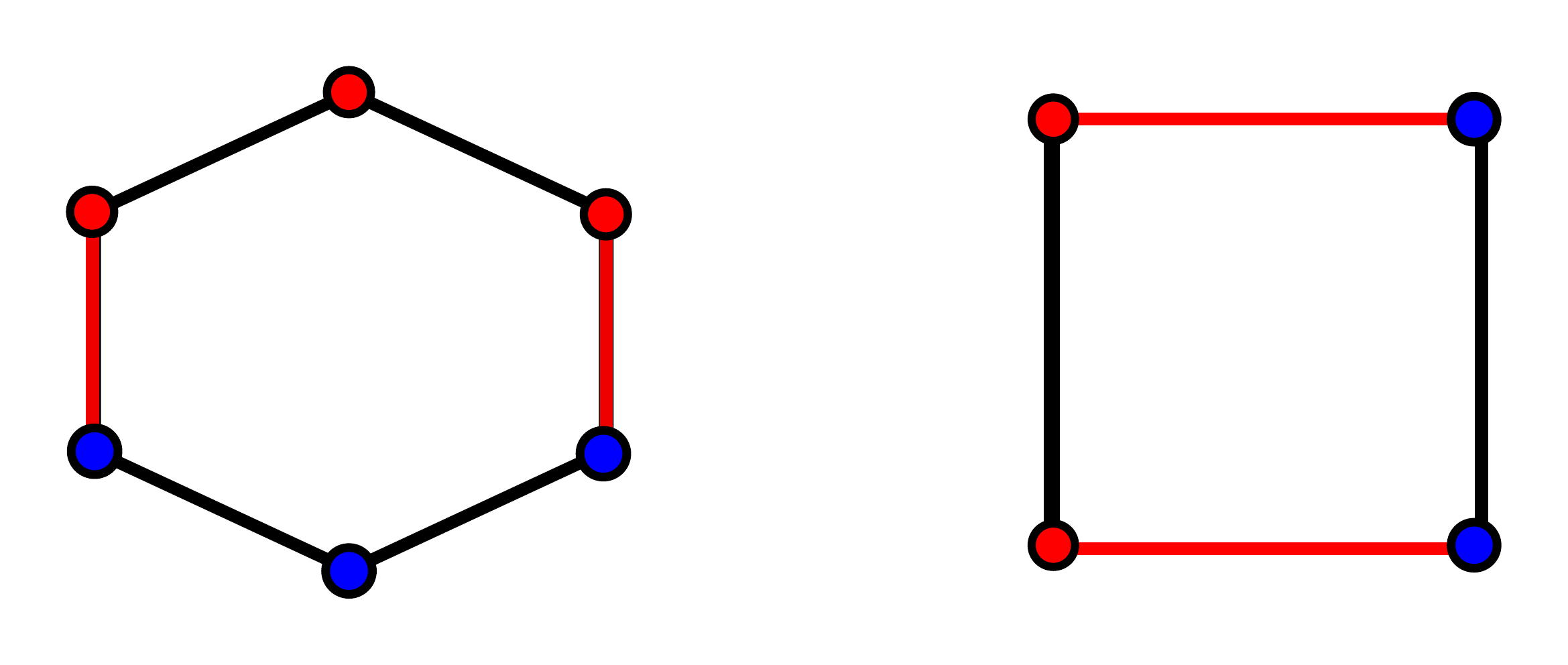}
    \caption{The 4AS-Square and 6AS-hexagon analyzed in Example \ref{subsec:A-Cells}}
    \label{fig:4CASAnd6ASPolygons}
\end{figure}

	Next, we analyze the 4AS-square.  This can be done completely without even drawing a handlebody diagram.  Here our situation is pictured Figure \ref{fig:4CASAnd6ASPolygons}, where vertices of the same color correspond to the same handlebody.  Suppose that $\Sigma$ has genus $g$.  Then we have the following:

	\begin{align*}
		\mathcal{X}^4(PQ) = \mathcal{X}^4(RS) = \natural^g(S^1 \times D^3) \\
		\mathcal{X}^4(QR) = \mathcal{X}^4(SP) = \natural^{g-1}(S^1 \times D^3) \\
		\mathcal{X}^4(PQRS) = \mathcal{X}^4(RSPQ) = \natural^{g-1}(S^1 \times D^3) \\
		\mathcal{X}^4_C(PQRSP) = \sharp^g(S^1 \times S^3) \\
		\mathcal{X}^4(QRSP) = \mathcal{X}^4(SPQR) = \natural^{g-1}(S^1 \times D^3) \natural(S^2 \times D^2) 
	\end{align*}
	The case of the 6AS-hexagon is completely analogous to the 4AS-square where we have the labeling in Figure \ref{fig:4CASAnd6ASPolygons}.

We now are in position to derive the following classical result from Hatcher's theorem that $\mathcal{P}(X)$ is simply-connected and our above analysis of the 2-cells:

\begin{theorem} \label{omega}
Every closed 4-manifold is cobordant to the connect sum of some number of $\mathbb{C}P^2$ and $\overline{\mathbb{C}P}^2$.
\end{theorem}

It then follows from the cobordism invariance of the signature that the oriented cobordism group $\Omega_4$ is isomorphic to $\mathbb{Z}$.  Note, as an aside, that the fact that $\mathcal{P}(\Sigma)$ is connected together with our construction of $\mathcal{M}^3(P_1,P_2)$ immediately yields that $\Omega_3 =0$.  

We will need to understand how $\mathcal{X}^4_C(L)$ changes when we alter $L$ by going over some 2-cell in $\mathcal{P}(\Sigma)$.  We begin with the following preliminary observation whose proof follows immediately from the definition of our construction:

\begin{lemma} \label{2-cell}
	Let $W_1,W_2$ be walks in $\mathcal{P}(\Sigma)$ with endpoints $P_1$ and $P_2$, and let $U$ and $V$ be two walks in $\mathcal{P}(\Sigma)$ so that the end point of $U$ is $P_1$ and the start point of $V$.  Then
	$$
	\mathcal{X}^4(UW_2V) = ( \mathcal{X}^4(UW_1V) - \mathcal{X}^4(W_1) ) \cup_{\mathcal{M}^3(P_1,P_2)} \mathcal{X}^4(W_2)
	$$
	This also holds in the case of $\mathcal{X}^4_C$ where the beginning of $U$ is the end of $V$.  
\end{lemma}

We will be applying Lemma \ref{2-cell}, in the case where $W_1 \cup W_2$ is the boundary of a 2-cell in $\mathcal{P}(\Sigma)$. This set up is pictured in Figure \ref{fig:HomotopingPath}.

Recall that a 1-surgery on a 4-manifold $X^4$ is the result of taking an embedding $\phi : S^1 \times D^3 \hookrightarrow X$ and forming the new 4-manifold
$$
X' = (X - \phi(\mathring{S^1 \times D^3})) \cup_{\partial \phi} (D^2 \times S^2)
$$
where we note that $\partial(D^2 \times S^2) = S^1 \times S^2 = \partial(S^1 \times D^3)$ and $\partial \phi$ denotes $\phi$ restricted to the boundary.  In this case we say that $X'$ is obtained from $X$ by a 1-surgery.  A 2-surgery on a 4-manifold is defined by switching the roles of $S^1 \times D^3$ and $D^2 \times S^2$ above.  Note that if $X'$ is obtained from $X$ by a 1-surgery, then $X'$ and $X$ are cobordant via the \emph{trace} of the surgery. Namely, given $\phi: S^1 \times D^3 \hookrightarrow X$, we can form
$$
(X \times I) \cup_{S^1 \times D^3 \subset X \times \{1\}} D^2 \times D^3
$$
which is a cobordism from $X$ to $X'$.  This can similarly be done if $X'$ is obtained from $X$ by a 2-surgery.  In our set up, we will not be seeing precisely the manifolds used in the definition of surgery, but nonetheless, the effect of cutting and pasting these pieces is simply a surgery. This is the content of the following lemma, where if we take $g=1$ in the following we get the usual definition of surgery.

	\begin{figure}
    \centering
   
    \labellist
    \pinlabel{n} at 670 80
    \pinlabel{g-1} at 250 120
    \endlabellist
     \includegraphics[scale=.3]{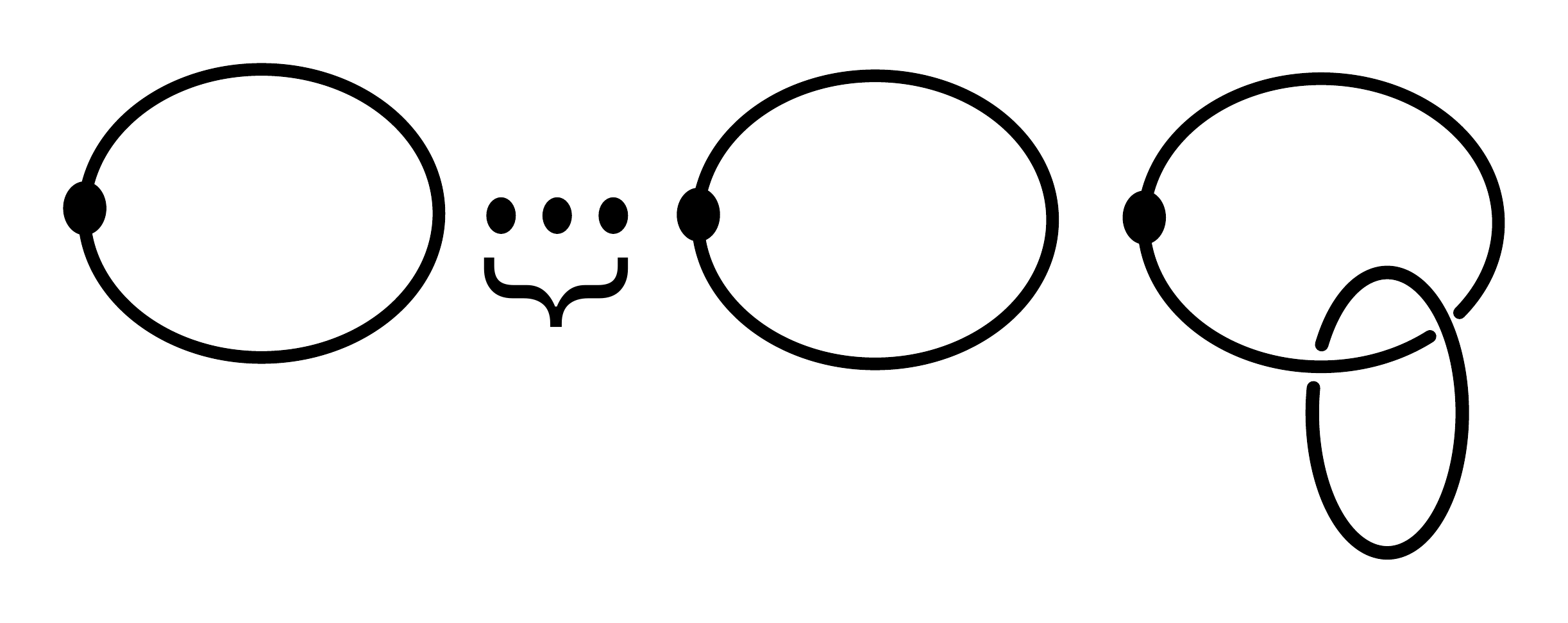}
    \caption{The handle decomposition of $\sharp^{g-1} (S^1 \times S^3)$ used for the surgery operation in Lemma \ref{surgery}}
    \label{fig:standardSurgery}
\end{figure}

\begin{lemma} \label{surgery}
	Let $X$ and $M$ be oriented 4-manifolds with $X \cong \sharp^{g-1} (S^1 \times S^3)$ given by the Kirby diagram in Figure \ref{fig:standardSurgery}. Let $X = X_1 \cup X_2$ where $X_1$ is the union of the 0- and 1-handles and $X_2$ is the union of the 2-, 3-, and 4-handles.  Let $\phi_1 : X_1 \hookrightarrow M$ and $\phi_2 : X_2 \hookrightarrow M$ be orientation-preserving inclusions.  Then
	$$
	(M - \phi_1(\mathring{X_1})) \cup_{\partial \phi_1} \overline{X_2}
	$$
and
$$
	(M - \phi_2(\mathring{X_2})) \cup_{\partial \phi_2} \overline{X_1}
$$
are obtained from $M$ by a single 1-surgery, and a single 2-surgery, respectively.  
\end{lemma}

\begin{proof}
We treat the case of 1-surgery explicitly, noting that the 2-surgery is just the inverse operation. Note that $X$ is constructed by first adding a cancelling 2-handle to $X_1$ followed by doubling what remains. When removing and reinserting the doubled portion, $M$ is unchanged. We are left with only removing a 1-handle, and replacing it by a 2-handle that cancels it, which is precisely the definition of surgery.
\end{proof}

If $X' \sharp \mathbb{C}P^2$ or $X' \sharp \overline{\mathbb{C}P}^2 = X$, then we call $X'$ a (+)- or (-)-blowdown of $X$, respectively. Likewise, we call $X$ a (+)- or (-)-blowup of $X'$. In the following construction, we will not be performing an operation corresponding to the definition of a blowups or blowdown. However, as in the previous lemma the overall effect is to perform blowups or blowdowns. The following lemma has a proof which is very similar to the previous lemma, so we omit it.

	\begin{figure}
    \centering
   
    \labellist
    \pinlabel{1} at 700 180
    \pinlabel{g-1} at 250 120
    \endlabellist
     \includegraphics[scale=.3]{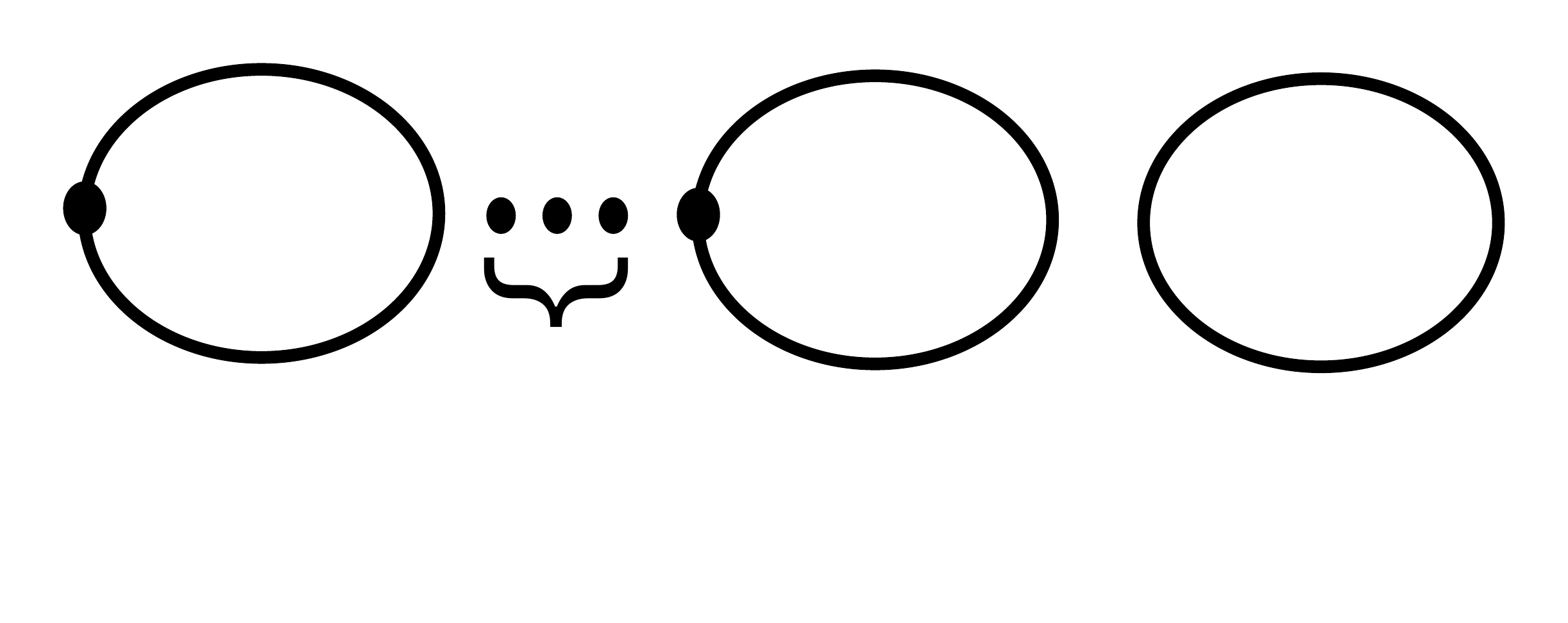}
    \caption{The handle decomposition of $\sharp^{g-1} (S^1 \times S^3) \sharp \mathbb{C}P^2$ used for the surgery operation in Lemma \ref{cp2}}
    \label{fig:standardCP2}
\end{figure}

\begin{lemma} \label{cp2}
	Let $X$ and $X'$ be oriented 4-manifolds with $X \cong \sharp^{g} (S^1 \times S^3) \sharp \mathbb{C}P^2$ given by the Kirby diagram in Figure \ref{fig:standardCP2}. Let $X = X_1 \cup X_2$ where $X_1$ is the union of the 0- and 1-handles and $X_2$ is the union of the 2-, 3-, and 4-handles.  Let $\phi_1 : X_1 \hookrightarrow X'$ and $\phi_2 : X_2 \hookrightarrow X'$ be orientation-preserving inclusions.  Then
	$$
	(X' - \phi_1(\mathring{X_1})) \cup_{\partial \phi_1} \overline{X_2}
	$$
and
$$
	(X' - \phi_2(\mathring{X_2})) \cup_{\partial \phi_2} \overline{X_1}
$$
	are a (-)-blowup, and a (-)-blowdown of $X'$, respectively.   
\end{lemma}

\begin{figure}
    \centering
    \includegraphics[scale=.25]{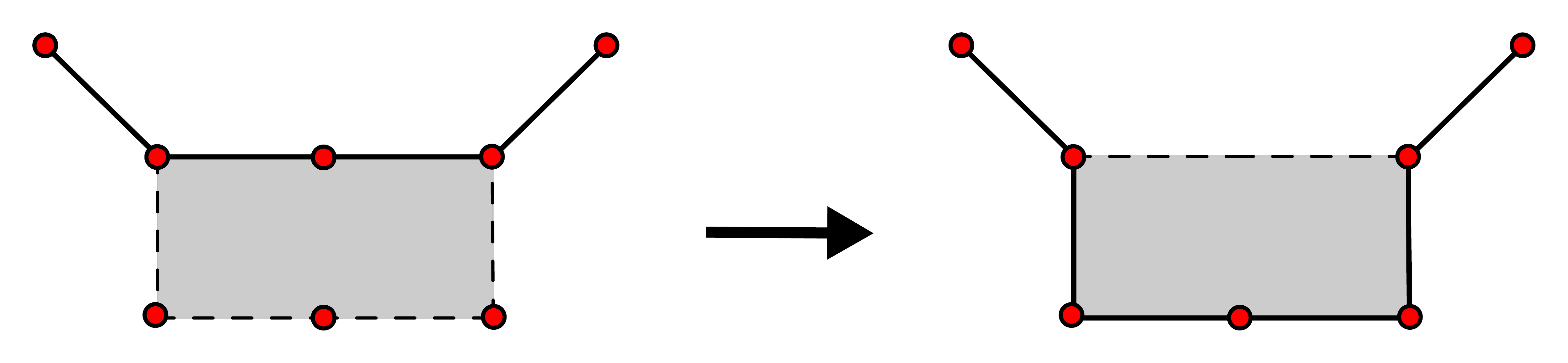}
    \caption{A homotopy of a path in a 2-complex consists of a sequence of operations of replacing a connected subset of the boundary of a 2-cell with the rest of the boundary.}
    \label{fig:HomotopingPath}
\end{figure}

Following the previous analysis we now understand how manifolds change as we homotope loops over 2-cells in the pants complex. This leads us to our main application.

\cobordismGroup

\begin{figure}
    \centering
    \labellist
		\pinlabel {The cobordism associated to} at 795 210
		\pinlabel {replacing one edge in a 3S} at 795 180
		\pinlabel {triangle by the other two edges} at 795 150
		
		\pinlabel {The cobordism associated to} at 780 575
		\pinlabel {replacing two edges in a} at 780 545
		\pinlabel {3S triangle by the other edge} at 780 515

		\pinlabel {The cobordism associated to} at -210 420
		\pinlabel {replacing portions of non-3S} at -210 385
		\pinlabel {2-cells by the other edge} at -210 355
		
		\pinlabel {We cap off the manifold } at -190 755
		\pinlabel {associated to the constant } at -190 730
		\pinlabel {path with $\natural{S^1 \times B^4}$} at -190 695
		
	\endlabellist
    \includegraphics[scale=.3]{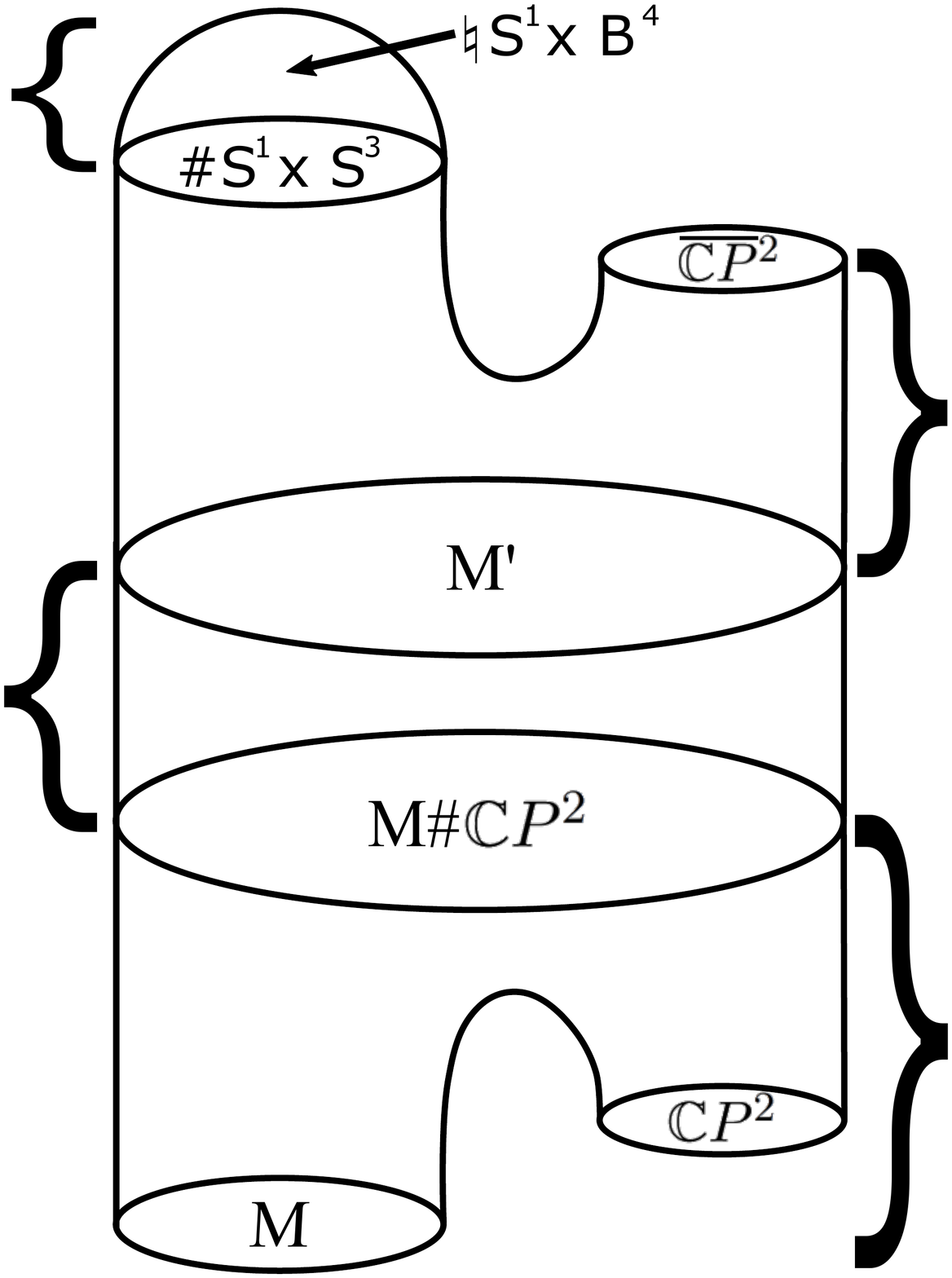}
    \caption{An illustration of the cobordism build in the proof of Theorem \ref{thm:cobordismGroup}}
    \label{fig:buildingCobordism}
\end{figure}

\begin{proof}

	Let $\Sigma$ be a genus $g$ surface, and let $L$ be a loop in $\mathcal{P}(\Sigma)$.  By Theorem \ref{thm:simplyConnected}, $\mathcal{P}(\Sigma)$ is simply-connected, so we know that there exists a cellular disk $D$ in $\mathcal{P}(\Sigma)$ with boundary $L$.  First, we alter $L$ using $D$ to a new loop called $L'$, so that $L'$ traverses a tree in $\mathcal{P}(\Sigma)$.  We can do this by changing the sides of the polygons that $L$ goes around, using $D$, as in Lemma \ref{2-cell}.    We will consider the different 2-cells following the order of the analysis in Section \ref{sec:2-cells}. 

	As analyzed in Subsection $\ref{subsec:3S}$, changing $L$ using a 3S-triangle results in the manifold changing by a (+)- or a (-)-blowup or blowdown, depending on the orientation of the triangle and the partition of the edges. We may achieve a cobordism from $X^4_C(M)$ to the resulting manifold by introducing a disjoint copy of  $\mathbb{C}P^2$ or $\overline{\mathbb{C}P}^2$ and forming the standard cobordism from a disjoint union to the connected sum, as illustrated in Figure \ref{fig:buildingCobordism}. 
	
	As analyzed in Subsection \ref{subsec:A-Cells}, changing $L$ using an 3A-triangle, an 5A-pentagon, or an 4A-square does not change the resulting 4-manifold at all. This modification of the loop therefore corresponds to the product cobordism on the manifolds.

	Suppose we are changing $L$ using a 4S-square.  If we are in the situation where the edges are partitioned into two sets of two adjacent edges, then the effect of the move is to remove some set 4-dimensional 1-handlebody, and then to reinsert the same 4-dimensional 1-handlebody, so by \cite{LP}, this does not change the resulting 4-manifold.  If we are in the case where the edges are partitioned into two sets where one set has three adjacent edges, then this affects $L$ by removing $\natural^{g-2}(S^1 \times D^3)$ and inserting $\natural^{g-2}(S^1 \times D^3) \natural (S^2 \times D^2)$ or vice versa, depending on which edges belong to $L$.  By Lemma \ref{surgery} this affects the resulting 4-manifold by performing either a 1- or a 2-surgery on $\mathcal{X}^4_C(L)$. We may achieve a cobordism between these manifolds using the trace of this surgery. The case of an 4AS-square is completely analogous.  

	After collapsing all of the 2-cells, we will have a loop $L'$ that is traversing a tree in $\mathcal{P}(\Sigma)$ with $\mathcal{X}^4_C(L')$, $\mathcal{X}^4_C(L)$, and some number of $\mathbb{C}P^2$'s and $\overline{\mathbb{C}P}^2$'s cobounding a 5-manifold.  We next collapse $L'$ to just a point.  To do this we proceed inductively on the size of the tree.  Choose a leaf of the tree, so that $L'$ must traverse this leaf in one direction, and then immediately turn back and go in the other direction.  Let $L''$ be the loop obtained from $L$ by removing this redundant edge followed by its reverse.  If the edge of the leaf is an A-edge, then $\mathcal{X}^4_C(L'')$ is equal to $\mathcal{X}^4_C(L')$.  So we need only consider the case where this leaf edge is an S-edge.

	Let $W$ be a walk in $\mathcal{P}(\Sigma)$ that is an S-edge traversed twice in a row in opposite directions.  The effect of removing the S-leaf from $L'$ is the same as removing a copy of $\mathcal{X}^4(W)$ from $\mathcal{X}^4_C(L')$ and replacing the result with $\natural^g(S^1 \times D^3)$, which is the manifold with boundary associated to the constant path.  But, in the same way that we analysed the boundaries of the 2-cells above, we find that $\mathcal{X}^4(W)$ is always $\natural^{g-1}(S^1 \times D^3) \natural (D^2 \times S^2)$.  But then, by Lemma \ref{surgery}, $\mathcal{X}^4_C(L'')$ is obtained from $\mathcal{X}^4_C(L')$ by performing a 1-surgery. We may again  achieve a cobordism between these manifolds using the trace of the surgery. Then by induction on the number of leaves, we may repeat this process until we arrive at the constant path. In the end, we find that there is a 5-manifold whose boundary is a disjoint union of $\mathcal{X}^4_C(L)$, some number of $\mathbb{C}P^2$ and $\overline{\mathbb{C}P^2}$, and $\sharp^g S^1 \times S^3$.  By capping off the $\sharp^g S^1 \times S^3$ with $\natural^g{S^1 \times B^4}$, the result follows.  

\end{proof}

\section{An explicit cobordism}

\begin{figure}
    \centering
    \includegraphics[scale=.6]{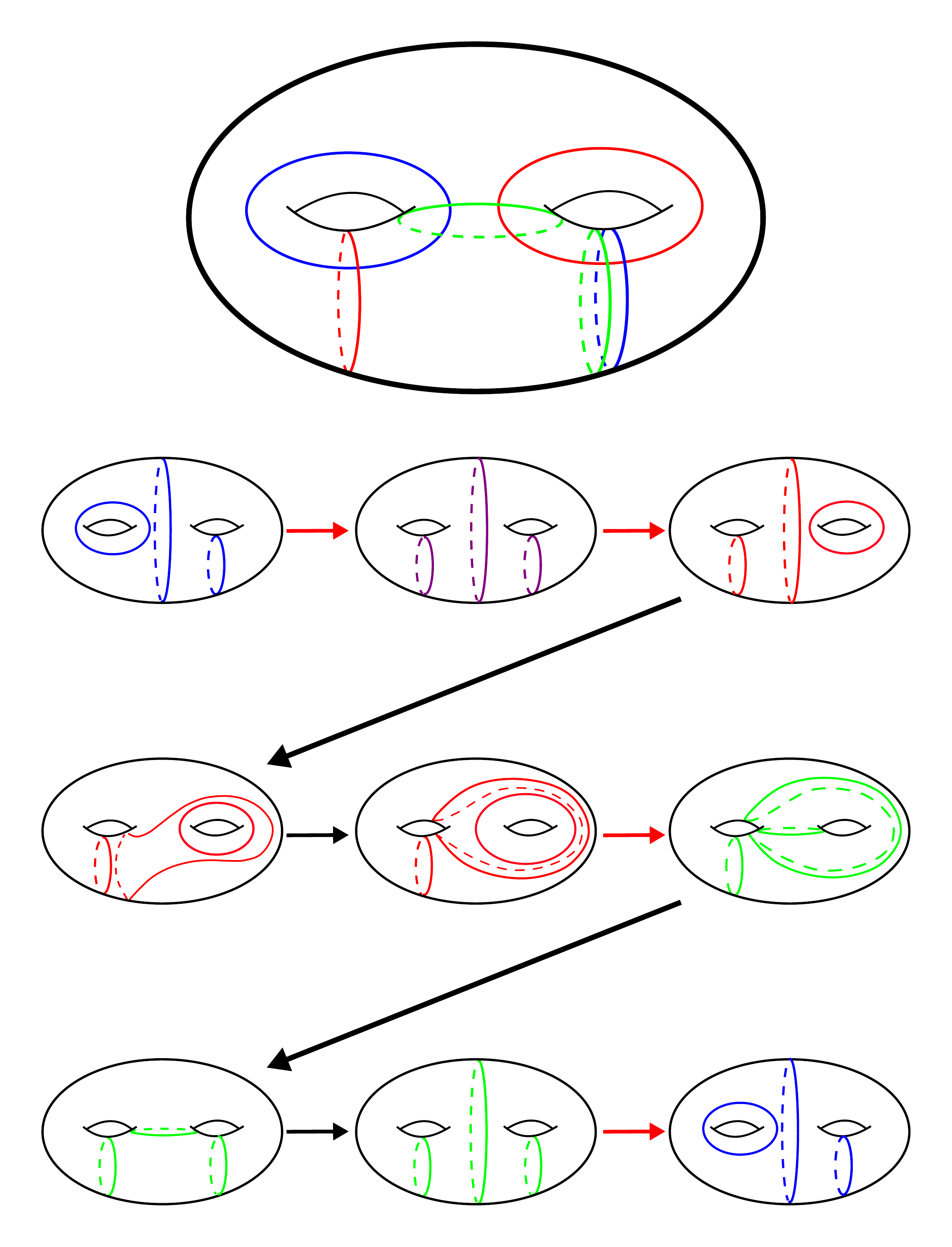}
    \caption{Top: An unbalanced trisection diagram for $S^4$. Bottom: A loop in the pants complex corresponding to this trisection.}
    \label{fig:S4PantsLoop}
\end{figure}

\begin{figure}
    \centering
    \includegraphics[scale=.5]{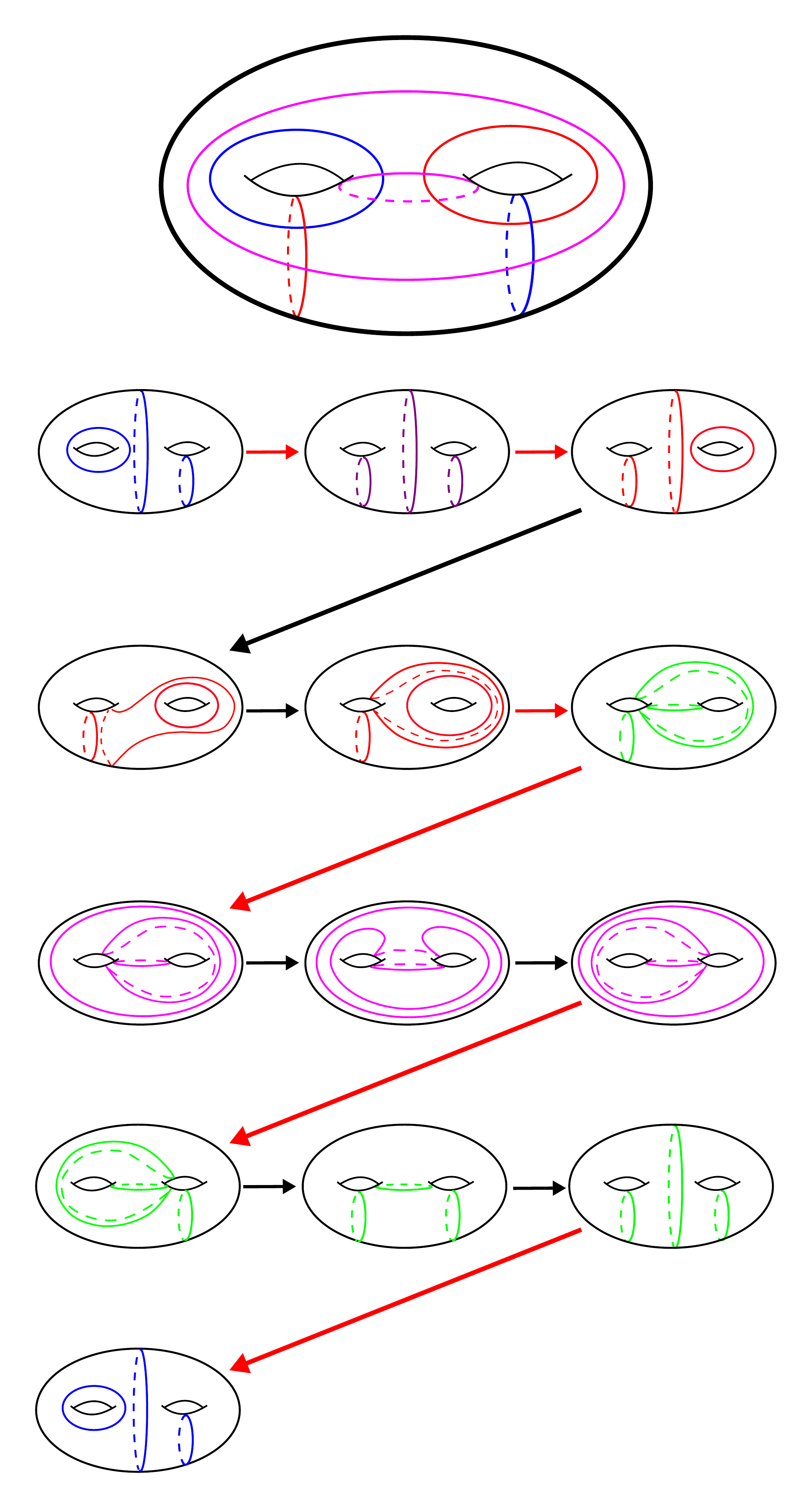}
    \caption{Top: A trisection diagram for $S^2 \times S^2$. Bottom: A loop in the pants complex corresponding to this trisection.}
    \label{fig:S2xS2PantsLoop}
\end{figure}

\begin{figure}

    \centering
       \labellist
		\pinlabel {\LARGE{\textbf{$6AS$}}} at 570 410
	\endlabellist
    \includegraphics[scale=.5]{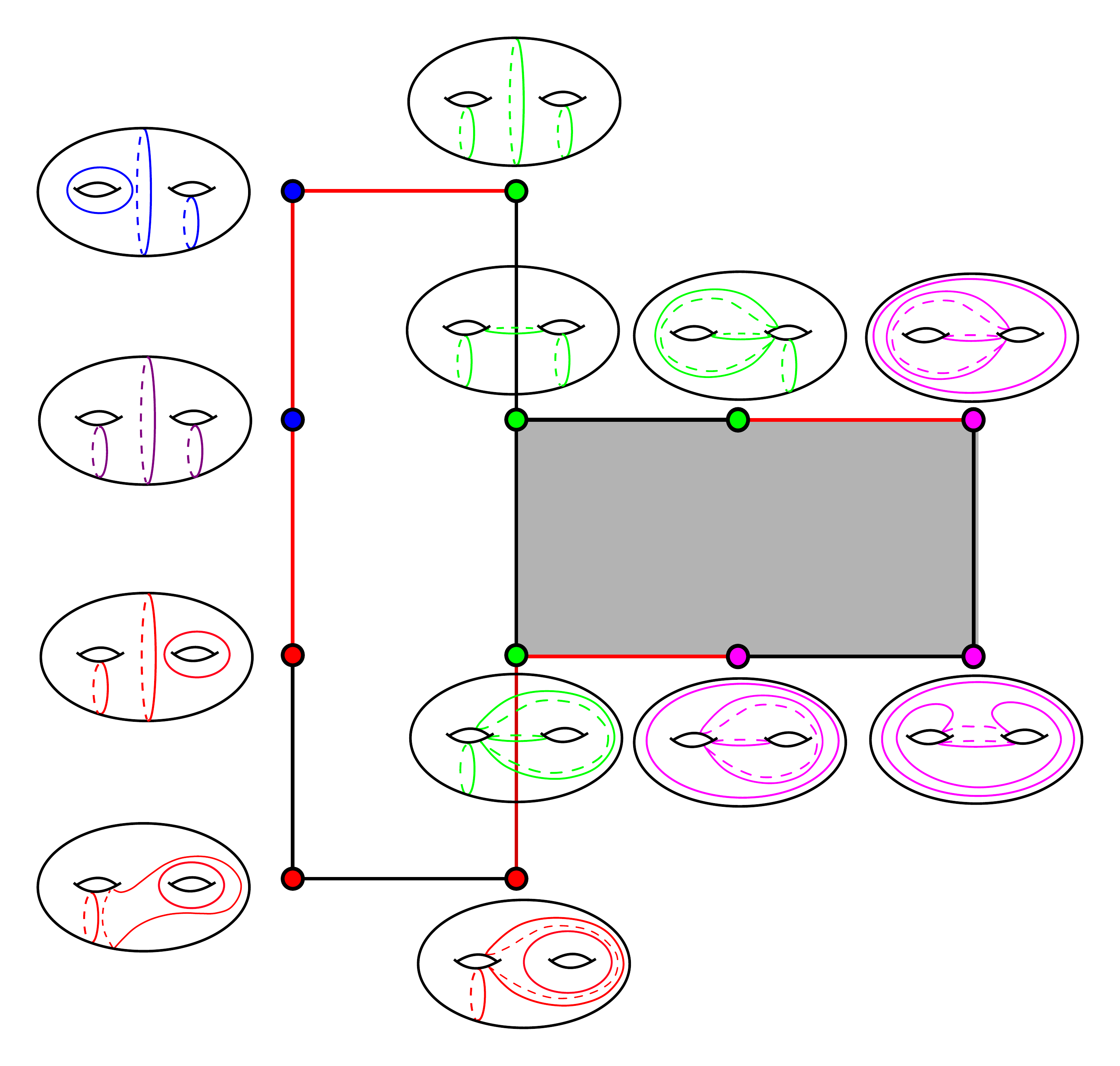}
 
    \caption{The loops for $S^2 \times S^2$ and $S^4$ meet in a $6AS$-hexagon. Taking the short way around the hexagon yields the loop for $S^4$ whereas the long way gives corresponds to $S^2 \times S^2$.}
    \label{fig:CobordismS4andS2xS2}
 
\end{figure}

In this section, we produce an explicit cobordism between $S^4$ and $S^2 \times S^2$ corresponding to a trivial 5-dimensional 2-handle attached to $S^4 \times [0,1]$ along $S^1 \times D^3 \subset S^4 \times \{1\}$. We start by describing loops in the pants complex corresponding to $S^4$ and $S^2 \times S^2$ coming from trisections of these manifolds. At the tops of Figures \ref{fig:S4PantsLoop} and \ref{fig:S2xS2PantsLoop} we see trisection diagrams for $S^4$ and $S^2 \times S^2$, respectively. These give rise to the loops in the pants complex shown below each of these diagrams. These loops meet in a 6AS hexagon as shown in Figure $\ref{fig:CobordismS4andS2xS2}$ and the two loops, together with the hexagon represent the desired cobordism.

\section{Signature}

In this section, we explain one way of computing the signature of $\mathcal{X}^4_C(L)$. By Novikov additivity, if $X$ and $X'$ are two 4-manifolds with diffeomorphic boundary, and $Y$ is a manifold obtained by gluing $X$ and $X'$ along their boundaries by a orientation-reversing diffeomorphism, then $\sigma(Y) = \sigma(X) + \sigma(X')$. When the manifolds $X_1$ and $X_2$ are not glued along their whole boundary, but rather some submanifolds of their boundary, then we no longer have this additivity. However, there is a correction term that was identified by Wall \cite{Wall}.  This correction term was further identified with the Maslov index of a certain triple of Lagrangians in \cite{cappell} (see also \cite{ranicki}).  It is in this form that we will apply it.

Let $(V,\psi)$ be a finite-dimensional vector space over $\mathbb{Q}$ together with a nonsingular symplectic form, and let $L_1,L_2,L_3 \subset V$ be three Lagrangians.  The Maslov index $M(L_1,L_2,L_3) \in \mathbb{Z}$ is the signature of the singular symmetric form given by
\begin{align*}
    \theta : L_1 \oplus L_2 \oplus L_3 \times L_1 \oplus L_2 \oplus L_3 &\to \mathbb{Q}\\
    ((x_1, x_2,x_3),(y_1,y_2,y_3)) &\mapsto \sum_{i \neq j}(-1)^{i+j}\psi(x_i,y_j)
\end{align*}

Suppose that the boundaries of $X$ and $X'$ have both been Heegaard split with the same genus surfaces so that $\partial X = H_1 \cup_\Sigma H_2$ and $\partial X' = H_1' \cup H_2'$.  Let $Y$ be the oriented 4-manifold that results from gluing $H_2$ to $H_2'$ by an orientation-reversing diffeomorphism $\phi$.  Let $L_1,L_2$, and $L_3$ be the Lagrangians in $H_1(\Sigma; \mathbb{Q})$ that are the kernels of the inclusions of $\Sigma$ as the boundary of $H_1, H_2$, and (using $\phi$) $H_2'$, respectively.  Then we have $\sigma(Y) = \sigma(X) + \sigma(X') - M(L_1,L_2,L_3)$.  

In our setting, if $L$ is an oriented loop in $\mathcal{P}(\Sigma)$, with vertices $P_1 P_2 \cdots P_n$, then applying the aforementioned formula, together with Novikov additivity when the last wedge is added and the fact that $\sigma(S^1 \times D^3) = 0$, we have the following:

\begin{proposition} Let $L$ be a loop in $\mathcal{P}(\Sigma)$, then
$$
\sigma(\mathcal{X}^4_C(L)) = -\sum_{i=2}^{n-1} M(L_1, L_i, L_{i+1})
$$
\end{proposition}

If $L$ is a loop in $P(\Sigma)$ we define $\sigma(L)$ to be the integer $\sigma({X}^4_C(L))$. We emphasize that the previous definition can be calculated using only information about the loop, with no reference to the 2-skeleton of the pants complex.

In this paper, we primarily have used information about the pants complex to derive information about 4-manifolds, but this process can also be reversed.  For example, given a loop $L$ in $\mathcal{P}(\Sigma)$, this loop bounds a disk $D$; what can we say about the 2-cells that make up $D$?  By orienting $L$, the disk $D$ inherits an orientation, and, in particular all of the 3S-triangles in $D$ inherit an orientation as well. Each triangle in this disk is either positive and negative, namely those that give rise to $\mathbb{C}P^2$ are positive, and those that give rise to $\overline{\mathbb{C}P}^2$ are negative. Following, Theorem \ref{thm:cobordismGroup}, we see that summing the number of positive 3S-triangles in $D$ and subtracting the negative 3S-triangles in $D$ is gives us $\sigma(L)$. We encapsulate the previous discussion in the following corollary.

\loopTriangleBound

\bibliography{thebib.bib}
\bibliographystyle{plain}

\end{document}